\documentclass[review]{elsarticle}

\usepackage[a4paper, total={6in,9in}]{geometry}

\usepackage{lipsum}
\usepackage{xcolor, hyperref}
\hypersetup{
    colorlinks = true,
    linkcolor = red,
    linkbordercolor = {white}
}

\usepackage{lipsum}
\usepackage[utf8]{inputenc}
\usepackage{amssymb}
\usepackage{amsmath}
\usepackage{graphicx} 
\usepackage{xspace}
\usepackage{multirow}
\usepackage{bbm}
\usepackage{algorithm}
\usepackage{algpseudocode}
\usepackage{graphicx}
\usepackage{subcaption}
\usepackage{enumitem}
\usepackage{multirow}
\usepackage{adjustbox}
\usepackage{bm}

\usepackage{lineno}
\modulolinenumbers[1]

\newcommand{\red}[1]{{#1}}

\begin{document}

\newcommand{\igor}[1]{\textcolor{red}{#1}}
\newcommand{\spann}{\text{span}}
\newcommand{\argmin}{\text{argmin}}

\begin{frontmatter}

\title{Stability estimates for radial basis function methods applied to linear scalar conservation laws}


\author[igor]{Igor Tominec}
\ead{igor.tominec@math.su.se}

\author[em]{Murtazo Nazarov}
\ead{murtazo.nazarov@it.uu.se}

\author[em]{Elisabeth Larsson}
\ead{elisabeth.larsson@it.uu.se}

\address[igor]{Stockholm University, Department of Mathematics, Division of Computational Mathematics}
\address[em]{Uppsala University, Department of Information Technology, Division of Scientific Computing}

\begin{abstract}
    We derive stability estimates for three commonly used radial basis function (RBF) methods to 
    solve hyperbolic time-dependent PDEs: the 
    RBF generated finite difference (RBF-FD) method, the RBF partition of unity method (RBF-PUM) and Kansa's (global) RBF method. 
    We give the estimates 
    in the discrete $\ell_2$-norm intrinsic to each of the three methods.
    The results show that Kansa's method and RBF-PUM can be $\ell_2$-stable in time under a sufficiently large oversampling 
    of the discretized system of equations. \red{The RBF-FD method in addition requires stabilization of the 
    spurious jump terms due to the discontinuous RBF-FD cardinal basis functions.}
    Numerical experiments show an agreement with our theoretical observations.
    \end{abstract}


\end{frontmatter}


\section{Introduction}
\label{Sec:intro}
This paper aims to explore the theory behind the stability properties of widely used radial basis function (RBF) methods applied to solve the following \red{linear} scalar conservation laws.  Find an unknown function $u:\mathbb R^2 \times \mathbb R^+ \mapsto \mathbb R^2$ such that:
\begin{equation}
      \label{eq:intro:conservationlaw}
      \begin{aligned}
\partial_t u &= - \nabla \cdot \bm f (u) \quad \mbox{in } \Omega \times (0, T],\\
u(\bm x, 0) &= u_0(\bm x) \quad \mbox{in } \Omega,
      \end{aligned}
\end{equation}
with appropriate boundary conditions, where $\Omega$ is an open and bounded domain in $\mathbb R^2$, $T>0$ is the final time, $\bm f \in \mathcal C^1(\mathbb R ; \mathbb R^{2})$ is the flux function, and $u_0(\bm x)$ is the initial data.

The theory of convergence analysis of methods used to solve \eqref{eq:intro:conservationlaw} requires that the solutions are stable in the $L_\infty$-norm, see e.g., \cite{DiPerna_1985}, which is usually referred to as maximum principle. The stability of the scheme with respect to the $L_2$-norm is a necessary step in proving $L_\infty$-estimates and the consistency of entropy inequalities, see for example the analysis presented in \cite{Nazarov13}. Therefore, in this paper, we will limit ourselves to analyzing and establishing the stability of some well-known RBF schemes with respect to the $L_2$-norm.

The three RBF methods that we investigate in this work are: $(i)$ the RBF-generated finite difference (RBF-FD) method, \cite{Tolstykh02, FBFN_book}; $(ii)$ the RBF partition of unity method (RBF-PUM) \cite{Wend02};  and $(iii)$ Kansa's RBF method \cite{KansaMethod}. The first two methods are localized such that the final system of equations is sparse, and the third method is global and gives a dense system of equations. Based on numerical experiments, it is known that all three methods are generally unstable when calculating collocated solutions of nonstationary conservation laws without additional stabilization, see e.g., \cite{FornbergLehto,Aiton14,Tominec17}. Stabilizations, mostly based on supplementing the numerical scheme by a hyperviscosity term, have been studied in \cite{Shankar_hypervi2, Shankar_hypervi1, FornbergLehto}. The authors of these references have successfully developed a hyperviscosity term by using higher powers of the Laplacian operator. They carefully constructed artificial hyperviscosity coefficients to stabilize linear advection-diffusion problems using both explicit and explicit-implicit time-stepping approaches. In their paper \cite{Shankar_semilagrangian}, the authors employ the semi-Lagrangian formulation of Kansa's method, RBF-PUM, and the RBF-FD methods to achieve stable numerical approximations for a linear transport problem over the surface of a sphere without utilizing any hyperviscosity terms.

Oversampling approximations have been shown to have better stability properties than the collocation counterpart for some RBF methods. For instance Kansa's method applied to \eqref{eq:intro:conservationlaw} becomes stable in time when the system of equations is oversampled instead of collocated, see \cite{RodrigoPlatte16}.
However, oversampled approximation does not make RBF-FD stable. In our recent work \cite{tominec21nonlconservation}, we showed that the RBF-FD method is generally unstable in time, both in the oversampled and in the collocated context. 

An alternative approach to achieving stable RBF approximations involves transforming them to the Galerkin formulation, i.e., writing the discretization in weak form. In the work \cite{Glaubitz1, Glaubitz2}, stability estimates were obtained for time-dependent linear conservation laws discretized using a Kansa-type method that uses a weak variational form through exact integration, and it was emphasized that the method requires imposing weak boundary conditions to be stable. Additionally, in their subsequent work \cite{Glaubitz3}, the authors presented summation-by-parts RBF operators adapted to the context of a Kansa-type method, which gives an energy stable discretization of the linear advection problem in one spatial dimension.

The main contribution of the present work is an investigation 
of the underlying reasons for instability of RBF discretizations when solving linear advection problems.  We find that 
(a) Kansa's method and RBF-PUM can be stable in time under a sufficient amount of oversampling, 
\red{(b) the RBF-FD method in addition requires stabilization of the spurious jump terms due to its discontinuous trial space (jumps in the cardinal basis functions),}
(c) the RBF-FD method can be stabilized by forcing the jump term over the interfaces of the discontinuities towards $0$,
(d) strong imposition of the inflow boundary conditions is in practice not preventing Kansa's method, RBF-PUM, or the stabilized RBF-FD method from being stable in time. 
We want to highlight that the jump stabilization from (c) requires the
definition of a dual mesh (Voronoi diagram) over the computational domain, since the discontinuities in the solution occur at the edges of that.

The paper is organized as follows: in Section \ref{sec:discretization} we discretize the nonlinear conservation law \eqref{eq:intro:conservationlaw} in an oversampled 
context, when the numerical solution is spanned by a set of global cardinal functions. 
Then we discuss how to construct the global cardinal basis function by three different methods: in Section \ref{sec:methods:rbfKansa} we formulate Kansa's method (i), in Section \ref{sec:methods:rbfpum} we formulate RBF-PUM (ii), 
in Section \ref{sec:methods:rbffd} we formulate the RBF-FD method (iii). Then in Section \ref{sec:stabilityproperties} 
we derive semi-discrete stability estimates for each method, 
provide the interpretation and propose an alternative, jump-based stabilization of the RBF-FD method. 
In Section \ref{sec:experiments} we validate the theoretical results by computing the eigenvalue spectra of the discretized advection operators for (i), (ii) and (iii) 
and show the numerical convergence results when solving \eqref{eq:intro:conservationlaw}, as well as the stability properties during a long-time simulation.
We make a final discussion in Section \ref{sec:finalremarks}.

\section{Collocation and oversampled discretization of a time-dependent linear advection problem}
\label{sec:discretization}
Here we discretize the linear advection problem \eqref{eq:intro:conservationlaw} in a collocated and an oversampled sense.

The domain $\Omega$ on which we solve \eqref{eq:intro:conservationlaw} is discretized using two point sets:
\begin{itemize}
    \item the nodal point set $X = \{x_i\}_{i=1}^N$ for generating the cardinal basis functions,
    \item the evaluation point set $Y = \{y_j\}_{j=1}^M$ for sampling the PDE \eqref{eq:intro:conservationlaw}, where $M=q N$.
\end{itemize}
That means: we have in total $N$ nodal points, and $q$ evaluation points are placed in the Voronoi region centered around each $x_i \in X$. Therefore, the total number of evaluation points are $M=qN$.

\begin{figure}[h!]
    \centering
    \begin{tabular}{ccc}
        \hspace{-0.3cm} \textbf{Point set $X$}  & \hspace{-0.3cm} \textbf{Perturbed point set $X$} & \hspace{-0.3cm} \textbf{Point sets $X$ and $Y$} \\
        \hspace{-0.3cm} \includegraphics[width=0.31\linewidth]{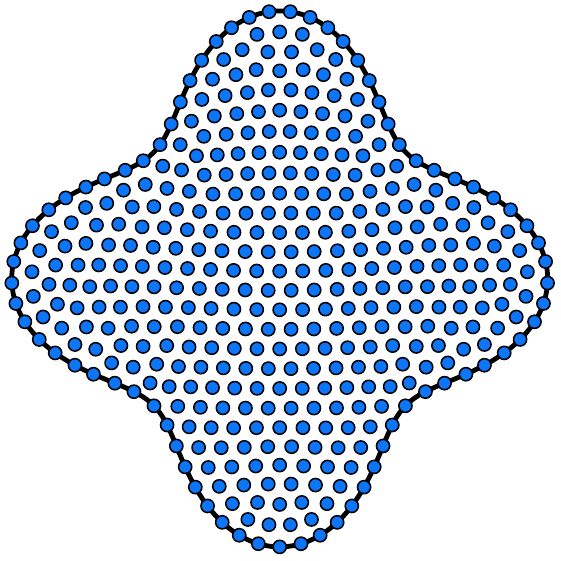}           &
        \hspace{-0.3cm} \includegraphics[width=0.31\linewidth]{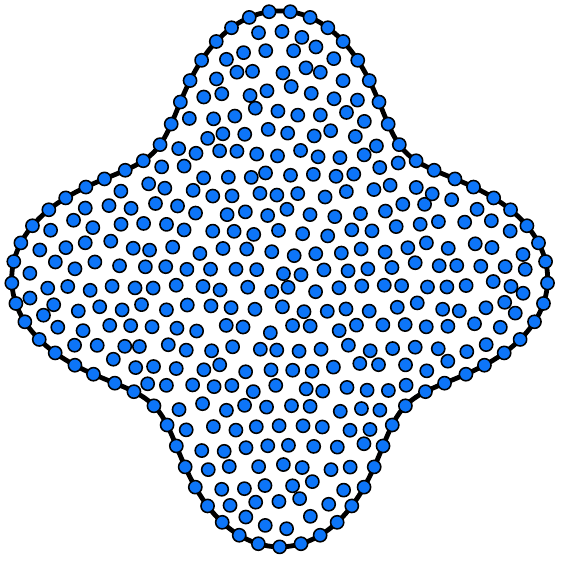} &
       \hspace{-0.3cm} \includegraphics[width=0.31\linewidth]{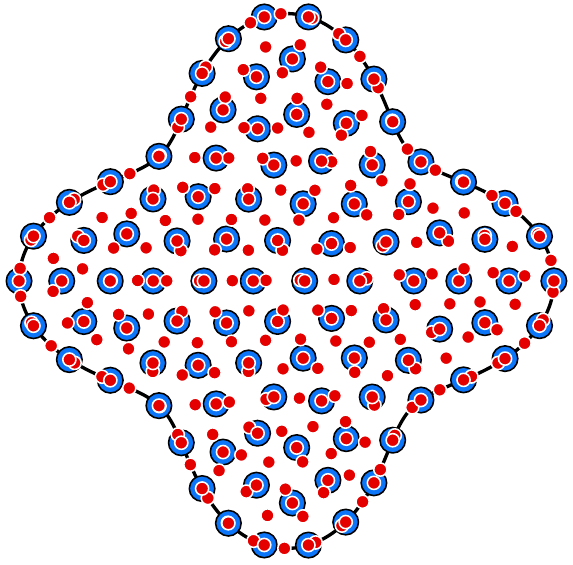}
    \end{tabular}
    \caption{Plots from left to right: a point set $X$ (blue markers) obtained using DistMesh, a randomly perturbed DistMesh $X$ point set (blue markers),
        point set $X$ (blue markers) together with the evaluation point set $Y$ (red markers), where the oversampling parameter is set to $q=4$.}
    \label{fig:experiments:pointset}
\end{figure}

An example of two point sets is visualized in Figure \ref{fig:experiments:pointset}.
In the present work, we generate the $X$ point set such that the mean distance between the points is set to $h$, by using the
DistMesh algorithm \cite{DistMesh}, which is in some cases randomly perturbed to also generalize the numerical results to 
less regular set of points. When using the oversampled discretization,
the $Y$ point set is also generated by DistMesh, but with an internodal distance $h_y = \frac{h}{\sqrt{q}}$. 
When using the collocation discretization we set $Y=X$, where, as a consequence we have that $h=h_y$.
We note that for a given point set $Y \subset \mathbb{R}^2$, 
we can define $h_y = ( \frac{|\Omega|}{M} )^{1/2}$. 
Here $|\Omega|$ is the area of $\Omega$, which can be estimated by $|\Omega| = \tilde M \frac{|\Omega_c|}{M_c}$, 
where $\Omega_c$ is a two dimensional box that encapsulates $\Omega$, $M_c$ is the number of points inside $\Omega_c$ and $\tilde M$ is 
the number of points in $\Omega_c \cap \Omega$.

To evaluate the RBF solution and its derivatives for any of the methods, we use:
\begin{equation}
    \label{eq:methods:ansatz}
    u_h(y, t) = \sum_{i=1}^N u_h(x_i,t) \Psi_i(y),\quad \mathcal{L}u_h(y, t) = \sum_{i=1}^N u_h(x_i,t)\, \mathcal{L}\Psi_i(y),
\end{equation}
where $\Psi_i(y)$, $i=1,..,N$, are the cardinal basis functions that together interpolate the unknown nodal values $u_h(x_i, t)$. The cardinal basis functions
are constructed differently and have different mathematical properties, depending on the choice of the RBF
method as shown in the sections below.

We next approximate the flux term as 
$
\bm f(u_h) = \sum_{i=1}^N \bm f(u_h(x_i,t)) \Psi_i(y_k)
$
for every $k=1,..,M$.
The discretization of \eqref{eq:intro:conservationlaw} is based on inserting the ansatz \eqref{eq:methods:ansatz} into
\eqref{eq:intro:conservationlaw}
and then sampling the PDE at the evaluation point set as follows:
\begin{equation}
    \label{eq:methods:advection_system}
    \begin{aligned}
        \sum_{i=1}^N \partial_t  u_h(x_i,t)\, \Psi_i(y_k) = - \sum_{i=1}^N \bm f(u_h(x_i, t)) \cdot \nabla \Psi_i(y_k), \\
        k=1,..,M.
    \end{aligned}
\end{equation}

Let us collect the nodal values of the flux function into the following matrix:
\[
\bm f(u_h(X,t)) := 
\left(
\begin{aligned}
 f_1(u_h(x_1,t)), f_1(u_h(x_2,t)), .., f_1(u_h(x_N,t))\\ 
 f_2(u_h(x_1,t)), f_2(u_h(x_2,t)), .., f_2(u_h(x_N,t)) 
 \end{aligned} 
 \right)^\top.
\]

Then, \eqref{eq:methods:advection_system} can be written on the matrix-vector format as:
\begin{equation}
    \label{eq:methods:advection_system_matrix_tmp}
    E_h \partial_t u_h(X,t) = \bm D_h \cdot \bm f(u_h(X,t)),
\end{equation}
where $(E_h)_{ki} = \Psi_i(y_k)$ and 
$(\bm D_h)_{ki} = -\nabla \Psi_i(y_k)$.
When $M=N$ in \eqref{eq:methods:advection_system}, then the discretization is performed in the collocation sense and the matrices $E_h$ and $\bm D_h$
in \eqref{eq:methods:advection_system_matrix_tmp} are square.
When $M>N$, the discretization is oversampled: $E_h$ and $\bm D_h$ are rectangular. In order to simplify the theoretical analysis at a later stage in this paper,
we introduce a norm scaling of \eqref{eq:methods:advection_system_matrix_tmp} for a 2D problem \cite{ToLaHe21},
and define $\bar E_h = h_y E_h$ and $\bar{\bm D}_h = h_y \bm D_h$, to obtain:
\begin{equation}
    \label{eq:methods:advection_system_matrix}
    \bar E_h \partial_t u_h(X,t) = \bar{\bm D}_h \cdot \bm f(u_h(X,t)).
\end{equation}

Given that the matrices in \eqref{eq:methods:advection_system_matrix} are typically rectangular, we require a projection onto a specific column space to obtain a square system. This is a necessary step for the resulting ODE system to evolve over time. We project \eqref{eq:methods:advection_system_matrix} onto the column space of $E_h$
by multiplying both sides of the PDE by $E_h^T$:
\begin{equation}
    \label{eq:methods:advection_system_matrix_projected}
    \begin{aligned}
        (\bar E_h^T \bar E_h)\, \partial_t u_h(X,t) & = (\bar E_h^T\, \bar{\bm D}_h) \cdot \bm f(u_h(X,t)),
    \end{aligned}
\end{equation}
The components of the matrix products involved in \eqref{eq:methods:advection_system_matrix_projected} are:
\begin{equation}
    \label{eq:methods:innerProducts1}
    (\bar E_h^T \bar E_h)_{ij} = (\Psi_i, \Psi_j)_{\ell_2},\quad (\bar E_h^T\, \bar{\bm D}_h)_{ij} = - (\Psi_i, \nabla \Psi_j)_{\ell_2},
\end{equation}
where $(u,v)_{\ell_2} = h_y^2 \sum_{k=1}^{M} u(y_k)\, v(y_k)$ is a discrete inner product. The $\ell_2$-inner product 
is an approximation of the $L_2$-inner product $(u,v)_{L_2} = \int_\Omega u v\, d\Omega$ \cite{ToLaHe21}. Interestingly,
if we had initially projected \eqref{eq:intro:conservationlaw} using the $L_2$-inner product,
instead of employing oversampling and a projection that led to the $\ell_2$-inner product,
the final discretization would have yielded a pure
Galerkin scheme.
Finally, the system of ODEs is obtained from \eqref{eq:methods:advection_system_matrix_projected} by inverting $\bar E_h^T \bar E_h$:
\begin{equation}
    \label{eq:methods:advection_system_matrix_projected_final}
\partial_t u_h(X,t)  = (\bar E_h^T \bar E_h)^{-1}\, (\bar E_h^T\, \bar{\bm D}_h) \cdot \bm f(u_h(X,t)),
\end{equation}
In practice, we do not invert $(\bar E_h^T \bar E_h)$ directly, but instead use an iterative solver for
\eqref{eq:methods:advection_system_matrix_projected} in each step of an explicit time-stepping algorithm. Since $(\bar E_h^T \bar E_h)$ is
symmetric and typically well conditioned, we use the conjugate gradient method (function \texttt{pcg()} in MATLAB).

In this work we consider the problem \eqref{eq:intro:conservationlaw} with a Dirichlet boundary condition specified
on the inflow part of $\partial \Omega$.
The boundary conditions in \eqref{eq:methods:advection_system_matrix_projected_final} are imposed in two ways: 
(i) when studying the eigenvalue spectra in the experimental part of this
paper, we impose a zero Dirichlet boundary condition by removing the components of the solution vector $u_h(X,t)$ 
associated with 
the solution values on the inflow boundary,
together with the corresponding columns and rows in the matrix $(\bar E_h^T \bar E_h)^{-1}\, (\bar E_h^T\, \bar{\bm D}_h)$, (ii) when computing the numerical solution in time, we impose the condition using the injection method, that is,
we overwrite the solution $u_h(X,t_i)$, $i=1,2,...$ with the corresponding inflow boundary value, after each 
step of the explicit time-stepping algorithm. The procedures (i) and (ii) are equivalent.

\section{Kansa's global RBF method}
\label{sec:methods:rbfKansa}
We work with a global interpolant defined on $\Omega$, where the interpolant is exact for radial basis functions $\phi_l(x)$, $l = 1,..,N$, and
the monomial basis, $\bar p_k(x)$, $k=1,..,m$. The interpolant is:
\begin{equation}
    \label{eq:methods:phspolybasis}
    \begin{aligned}
        u_h(x,t) & =\sum_{l=1}^N c_l(t) \phi_l(x) + \sum_{k=1}^{m} \beta_k(t) \bar p_k(x),  \\
                 & \text{ subject to  } \sum_{l=1}^N c_l(t) \bar p_k(x_l)=0,\quad k=1,..,m,
    \end{aligned}
\end{equation}
where $c_l(t)$ are the unknown interpolation coefficients, $\beta_k(t)$ the unknown Lagrange multipliers, and where the number of monomial terms is
$m = \binom{D_m+2}{2}$, where $D_m$ is the degree of the monomial basis. The discussion in this paper applies for any RBF, however, in our computations we use
the cubic polyharmonic spline (PHS) basis $\phi_l(x) = \| x - x_l\|_2^3$. 
\red{PHS and the monomial basis as a combination is well described in \cite{03_Buhmann_book,07_Fasshauer_book,05_Wendland_book}. 
In the context of the RBF-FD method the basis is also studied in \cite{Barnett15,Bayona19,BFFB17}, where the focus is on high-order polynomial approximation.}
In \eqref{eq:methods:phspolybasis} the objective is now to determine the interpolation coefficients $c_l$ and the Lagrange multipliers $\beta_k$,
by requiring the interpolation conditions $u_h(x_i,t) = u(x_i,t)$, $i=1,..,N$, where $x_i \in X$ and $u(x_i,t)$ are the unknown nodal values 
to be fulfilled. The result is
a system of equations in matrix-vector format:
\begin{equation}
    \label{eq:methods:M}
    \underbrace{\begin{pmatrix}
            A   & P \\
            P^T & 0
        \end{pmatrix}}_{:=\tilde A}
    \begin{pmatrix}
        \underline{c}(t) \\
        \underline{\beta}(t)
    \end{pmatrix}
    =
    \begin{pmatrix}
        \underline{u}(X,t) \\
        0
    \end{pmatrix}\quad\Leftrightarrow\quad
    \begin{pmatrix}
        \underline{c}(t) \\
        \underline{\beta}(t)
    \end{pmatrix}
    =
    \tilde A^{-1}
    \begin{pmatrix}
        \underline{u}(X,t) \\
        0
    \end{pmatrix}.
\end{equation}
Here $A_{jl} = \phi_l(x_j)$ for indices $j,l=1,..,N,$ and $P_{jk} = \bar p_k(x_j)$ for indices $k=1,..,m$.
The solution $u_h(y,t)$, $y\in \Omega$, is then written by reusing the computed coefficients inside \eqref{eq:methods:phspolybasis} 
and only keeping the first $N$ terms, and disregarding the $N+1,..,N+m$ terms which are related to the Lagrange multipliers:
\begin{equation}
    \label{eq:method:Kansa_cardinalF}
    \begin{aligned}
        u_h(y,t) & =
        \underbrace{
            \begin{pmatrix}
                \phi_1(y),..,\phi_N(y),\,  \bar p_1(y),..,\bar p_m(y)
            \end{pmatrix}}_{:=b(y)}
        \begin{pmatrix}
            \underline{c}(t) \\
            \underline{\beta}(t)
        \end{pmatrix}
        =\left(b(y) \tilde A^{-1} \right)_{1:N}\, \underline{u}(X,t) =                                                          \\
                 & = \sum_{k=1}^{N} u(x_k,t) \left[ b(y,X)\, \tilde A^{-1} (X, X) \right ]_k \equiv \sum_{k=1}^N u(x_k,t) \Psi_k (y).
    \end{aligned}
\end{equation}
Here we also derived a vector, continuous in $y$, of cardinal basis functions $\Psi_1(y),..,\Psi_N(y)$, used to construct $u_h(y,t)$ in \eqref{eq:methods:ansatz}.
For Kansa's method, the matrices of the discretization
\eqref{eq:methods:advection_system_matrix} are dense. 
\begin{figure}
    \begin{tabular}{ccc}
        \hspace{-0.25cm}\textbf{Kansa's method} & \hspace{-0.25cm}\textbf{RBF-PUM} & \hspace{-0.25cm}\textbf{RBF-FD} \\
        \hspace{-0.25cm}\includegraphics[width=0.32\linewidth]{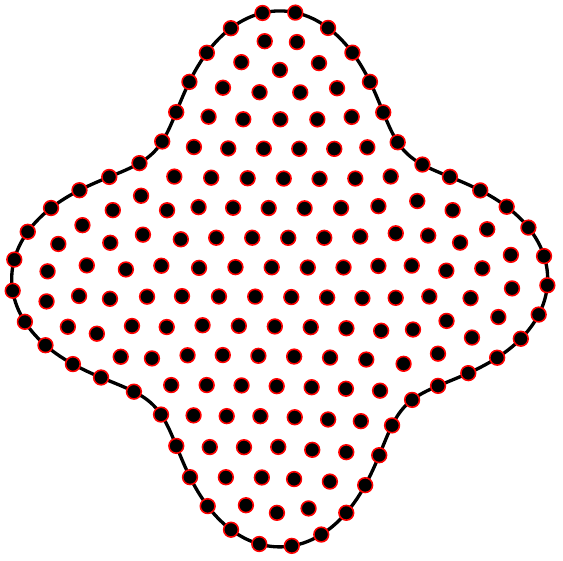} &
    \hspace{-0.25cm}\includegraphics[width=0.32\linewidth]{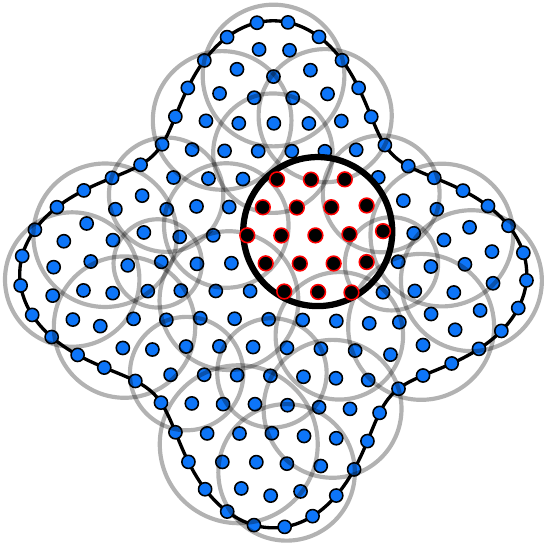} &
    \hspace{-0.25cm}\includegraphics[width=0.32\linewidth]{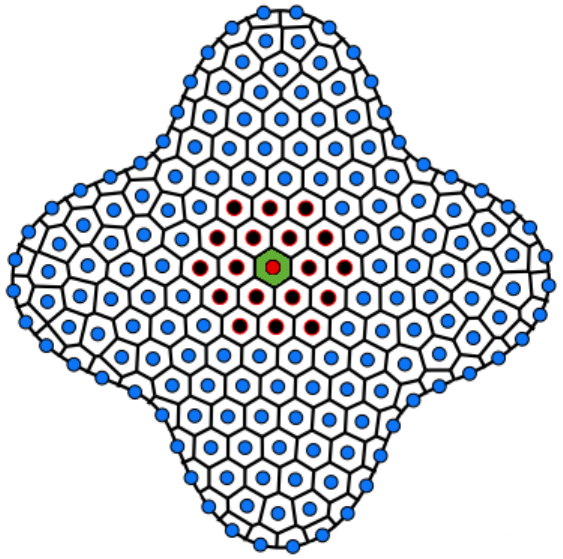}
    \end{tabular}
\caption{Nodes (blue points) and different node supports (black points with red edges) specific to three RBF methods: Kansa's RBF method (global support), 
the RBF partition of unity method (patch support) and the RBF generated finite difference (RBF-FD) method (stencil support). The green area in the RBF-FD method case 
illustrates the region to which the stencil approximation is further restricted when evaluating the numerical solution.}
\end{figure}

 \section{The RBF partition of unity method (RBF-PUM)}
\label{sec:methods:rbfpum}
As for Kansa's method, the domain $\Omega$ is first discretized using the nodal points $X$. 
Then we construct an open cover of $\Omega$
using the patches $\Omega^{(j)}$, $j=1,..,N_p$, such that $\Omega \subset \cup_{i=1}^{N_p} \Omega^{(j)}$. 
Each patch contains $n^{(j)} \geq 2 \binom{D_m+2}{2}$ points, where the points of the $j$-th patch are denoted by 
$X^{(j)} \subset X$. The RBF-PUM solution is represented as follows. We form an interpolant on $\Omega^{(j)}$ analogously to
Kansa's method on the whole $\Omega$: we use \eqref{eq:methods:ansatz} that leads to \eqref{eq:method:Kansa_cardinalF}, where $N$ is replaced by $n^{(j)}$,
and $X$ by $X^{(j)}$. The result is a local solution:
\begin{equation}
    \label{eq:methods:rbfpum_localsolution}
    u_h^{(j)}(y,t) = \sum_{k=1}^{n^{(j)}} u(x_k^{(j)},t) \left[ b(y,X^{(j)})\, \tilde A^{-1} (X^{(j)}, X^{(j)}) \right ]_k \equiv \sum_{k=1}^{n^{(j)}} u(x_k^{(j)},t) \psi_k^{(j)} (y),
\end{equation}
where $\psi_k^{(j)}(y)$ are now the local RBF-PUM cardinal functions. Then the global solution $u_h(y,t)$ is defined
by blending $u_h^{(j)}(y,t)$ together using compactly supported weight functions $w_j(y)$ such that:
\begin{equation}
    \label{eq:methods:rbfpum_globalsolution}
    u_h(y,t) = \sum_{j=1}^{N_p} w_j(y) u_h^{(j)}(y,t) = \sum_{j=1}^{N_p} w_j(y) \sum_{k=1}^{n^{(j)}} u(x_k^{(j)},t) \psi_k^{(j)} (y),
\end{equation}
where each $w_j(y)$ is constructed using Shepard's method \cite{Shepard}:
\begin{equation}
    w_j(y) = \frac{\Phi_j(y)}{\sum_{i=1}^{N_p} \Phi_i(y)},\quad \Phi_j(y)=(4r_j + 1)(1-r_j^4)_+,\quad r_j=\frac{\|y-\xi_j\|_2}{R_j}, \quad j=1,..,N_p.
\end{equation}
Here $\Phi_j$ are chosen as $\mathcal C^2(\Omega)$ Wendland functions \cite{WendlandFunctions} with compact support on $\Omega^{(j)}$, 
and $\xi_j$ and $R_j>0$ are consecutively the center and the radius of $\Omega^{(j)}$.
The global solution \eqref{eq:methods:rbfpum_globalsolution} is now cast in an equivalent form that is analogous to \eqref{eq:methods:ansatz}, by
introducing: (a) a set $J(i)$ of all patches where $x_i$ is contained, (b) an index operator $\kappa(j,i)$, 
which gives the local index of $x_i$ inside a patch $\Omega^{(j)}$.
We then have:
\begin{equation}
    \begin{aligned}
        \label{eq:methods:rbfpum_globalsolution2}
        u_h(y,t) & = \sum_{j=1}^{N_p} w_j(y) \sum_{k=1}^{n^{(j)}} u(x_k^{(j)},t) \psi_k^{(j)} (y)                                                                        \\
                 & = \sum_{i=1}^N u_h(x_i,t) \sum_{j \in J(i)} w_j(y) \psi_{\kappa(j,i)}^{(j)}(y) \equiv \sum_{i=1}^N u_h(x_i,t) \Psi_i(y).
    \end{aligned}
\end{equation}
and thus $\Psi_i(y) = \sum_{j \in J(i)} w_j(y) \psi_{\kappa(j,i)}^{(j)}(y)$.
By using the RBF-PUM cardinal basis functions $\Psi_i(y)$ to discretize a time-dependent PDE as in \eqref{eq:methods:advection_system},
the involved matrices in \eqref{eq:methods:advection_system_matrix} become sparse.

\section{The RBF generated finite difference (RBF-FD) method}
\label{sec:methods:rbffd}
The domain $\Omega$ is discretized using the nodal points $X$. 
Then we form a stencil for each $x_i \in X$, based on the $n = 2 \binom{D_m+2}{2}$ nearest neigbors, 
where each stencil is represented with a point set $X^{(j)} \subset X$, $j=1,..,N$.
On each stencil point set $X^{(j)}$ we form a local interpolation problem by using \eqref{eq:methods:ansatz} that leads to \eqref{eq:method:Kansa_cardinalF},
where $N$ is replaced by $n$,
and $X$ by $X^{(j)}$. The result is a set of stencil based approximations:
\begin{equation}
    \label{eq:methods:rbffd_localsolution}
    \begin{aligned}
        u_h^{(j)}(y,t) = \sum_{k=1}^{n} u(x_k^{(j)},t) \left[ b(y,X^{(j)})\, \tilde A^{-1} (X^{(j)}, X^{(j)}) \right ]_k \equiv \sum_{k=1}^{n} u(x_k^{(j)},t) \psi_k^{(j)} (y), \\
        j=1,..,N,
    \end{aligned}
\end{equation}
where $\psi_k^{(j)}(y)$ are the local RBF-FD cardinal functions. To arrive at the global solution $u_h(y,t)$, where $y \in \Omega$, we first
associate every $y$ with the index of the closest stencil center point defined as:
\begin{equation}
    \rho(y) = \arg \min_{i \in [1,N]} \| y- x_i \|_2.
\end{equation}
Analogously to the situation in RBF-PUM \eqref{eq:methods:rbfpum_globalsolution} we write the global solution as:
\begin{equation}
    \label{eq:methods:rbffd_globalsolution}
    u_h(y,t) = u_h^{\left(\rho(y)\right )}(y,t) = \sum_{k=1}^{n} u(x_k^{\left(\rho(y)\right )},t) \psi_k^{\left(\rho(y)\right )} (y).
\end{equation}
This implies that $u_h(y,t)$ restricted to a Voronoi region centered around $x_{\rho(y)} \in X$
is represented by a local solution $u_h^{\rho(y)}(y,t)$. 
An equivalent form of \eqref{eq:methods:rbffd_globalsolution} that is analogous to \eqref{eq:methods:ansatz} 
is written by using an index operator $\kappa(\rho(y), i)$ which in the RBF-FD method case gives the local index of $x_i$ 
inside a stencil point set $X^{(\rho(y))}$.
\begin{equation}
    \begin{aligned}
        \label{eq:methods:rbffd_globalsolution2}
        u_h(y,t) & = \sum_{k=1}^{n}
        u(x_k^{\left(\rho(y)\right )},t)\, \psi_k^{\left(\rho(y)\right )} (y)                                                                        \\
                 & = \sum_{i=1}^N u(x_i,t)\, \psi^{(\rho(y))}_{\kappa(\rho(y), i)} (y) \equiv \sum_{i=1}^N u(x_i,t) \Psi_i(y).
    \end{aligned}
\end{equation}
The global RBF-FD cardinal functions are thus $\Psi_i(y) = \psi^{(\rho(y))}_{\kappa(\rho(y), i)} (y)$. 
In this case we do not have any blending weight functions that ensure the continuity across the Voronoi regions, as is the case for RBF-PUM 
(see \eqref{eq:methods:rbfpum_globalsolution2}), 
and thus the RBF-FD global solution is discontinuous across the edges between these regions. 
\cite[Figure 2]{ToLaHe21}.

A MATLAB code for obtaining the collocated or the oversampled RBF-FD matrices is available from \cite{Tominec_rbffdcode2021}.

\section{Stability properties of RBF methods for linear advection problems}
\label{sec:stabilityproperties}
Let us consider the case when the flux is advective, i.e., $\bm f(u(x,t)) := \bm \beta(x,t)\, u(x,t)$, where $\bm \beta(x,t)$ is a vector field in $\mathbb R^2$. To simplify the below discussion we assume that this vector field is divergence-free: $\nabla \cdot \bm \beta \equiv 0$ for any $x\in \Omega$. Then, we can rewrite \eqref{eq:intro:conservationlaw} as
\begin{equation}\label{eq:adv}
      \begin{aligned}
		\partial_t u &= - \bm \beta \cdot \nabla u \quad \mbox{in } \Omega \times (0, T],\\
		u &= g \quad \mbox{in } \Gamma_{\text{inflow}} \times (0, T],\\
		u(\bm x, 0) &= u_0(\bm x) \quad \mbox{in } \Omega,
      \end{aligned}
\end{equation}
where $g(x,t)$ is a given sufficiently smooth function. 

We discretize \eqref{eq:adv} and following the steps discussed in Section~\ref{sec:discretization}, we obtain the following system of ODEs:
\begin{equation}\label{eq:adv:ode}
\bar E_h^T \bar E_h \partial_t u_h(X,t)  = (\bar E_h^T\, \bar{\bm D}_h) \cdot \bm \beta\, u_h(X,t).
\end{equation}

In this section, we derive semi-discrete $\ell_2$-stability estimates for the PDE discretization given in \eqref{eq:adv},
where the cardinal basis functions $\Psi_i$, $i=1,..,N$, are constructed using three RBF methods: Kansa's method described in Section \ref{sec:methods:rbfKansa},
RBF-PUM described in Section \ref{sec:methods:rbfpum} and the RBF-FD method described in Section \ref{sec:methods:rbffd}.

\subsection{Definitions of the inner products and norms}
\label{sec:stabilityproperties:innerproducts}
To measure the stability, we define a discrete inner product for functions $u_h$, $v_h$ sampled at $Y$-points:
\begin{equation}
    \label{eq:stabilityestimate:discreteinnerproducts}
    \begin{aligned}
        (u_h(Y,t), v_h(Y,t))_{\ell_2(\Omega)} & = (\bar E u_h(X,t), \bar E\, v_h(X,t))_{\ell_2(\Omega)}                           \\
                                          & = u_h^T(X,t) \bar E^T \bar E\, v_h(X,t) = h_y^2 \sum_{j=1}^M u_h(y_j,t)\, v_h(y_j,t),
    \end{aligned}
\end{equation}
where $\bar E$ is defined in the scope of \eqref{eq:methods:advection_system_matrix} and $h_y$ in Section \ref{sec:discretization} (we note that as $M\to \infty$, then $h_y \to 0$). 
The induced discrete norm is then:
\begin{equation}
    \label{eq:stabilityestimate:discretenorms}
    (u_h(Y,t), u_h(Y,t))_{\ell_2(\Omega)} =  u_h^T(X,t) \bar E^T \bar E\, u_h(X,t) = \| u_h(Y,t) \|_{\ell_2(\Omega)}^2.
\end{equation}
A link between our discretization and the chosen inner product is described in Section \ref{sec:discretization}.
The stability estimates derived later in this section are based on a relation between the $\ell_2$-inner product 
and the $L_2$-inner product. The latter is defined by:
\begin{equation}
    (u_h(y,t), v_h(y,t))_{L_2(\Omega)} = \int_\Omega u_h(y,t) v_h(y,t)\, d\Omega, \quad y\in \Omega,
\end{equation}
Let a Voronoi region $K_i$ be defined as:
\begin{equation}
    \label{eq:voronoi_region_Ki}
    K_i = \{ y\in\Omega\, |\, \|y-x_i\| \leq \|y-x_j\|,\, j \neq i,\, j=1,..,N\},\quad i=1,..,N.
\end{equation}
where we also note that $\cup_{i=1}^N K_i = \Omega$.
Now we use the integration error estimate \eqref{appendix:integral:Omega} derived in \ref{sec:appendix:integration}, 
to relate the two inner products as:
\begin{equation}
    \label{eq:stabilityestimate:discreteToContinuousInnerProduct}
    \begin{aligned}
    (u_h, v_h)_{\ell_2(\Omega)} &\leq \underbrace{C_{\int}\, h^{-1} h_y \max_i \|u_h\,v_h\|_{L_\infty(K_i)}}_{=e_I(u_h,v_h)} +  (u_h, v_h)_{L_2(\Omega)}    \end{aligned}
\end{equation}
where $C_{\int}$ is independent of $h_y$ and $h$. We have that for a fixed $h$, the integration error $e_I(u_h,v_h)$ approches $0$ as $h_y \to 0$.
The integration error estimate is later also numerically validated in Section \ref{sec:experiments:integration_error}.

\subsection{A semi-discrete stability estimate for Kansa's RBF method}
\label{sec:stabilityproperties:Kansa}
Consider a numerical solution $u_h$ spanned by the cardinal basis functions defined in \eqref{eq:method:Kansa_cardinalF}.
Observe that the regularity of $u_h$ is given by the lowest regularity of the bases $\phi_l$ and $\bar p_k$ used in \eqref{eq:method:Kansa_cardinalF}.
This implies $\Psi_i \in \mathcal C^2(\Omega)$ and $u_h \in \mathcal C^1([0,T]) \times \mathcal C^2(\Omega)$.
The derivation that follows is based on that regularity fact, since we are then allowed to differentiate the solution $u_h$ over the whole $\Omega$.

We consider the discrete problem \eqref{eq:methods:advection_system_matrix_projected} and multiply it with $u_h^T(X,t)$:
\begin{equation}
    \label{eq:stabilityKansa:firstStep}
    \begin{aligned}
        u_h^T(X,t) \bar E_h^T\, \bar E_h\, \partial_t u_h(X,t) & = u_h^T (X,t) \bar E_h^T\, \bar D_h\, u_h(X,t).
    \end{aligned}
\end{equation}
Using the relation \eqref{eq:stabilityestimate:discreteinnerproducts} for both sides of \eqref{eq:stabilityKansa:firstStep}
and using \eqref{eq:methods:innerProducts1} for the right-hand-side of \eqref{eq:stabilityKansa:firstStep}, we rewrite \eqref{eq:stabilityKansa:firstStep}
on the form:
\begin{equation}
    \label{eq:stabilityKansa:secondStep}
    \begin{aligned}
        (u_h(Y,t), \partial_t\, u_h(Y,t))_{\ell_2(\Omega)} & = - (u_h(Y,t), \bm \beta \cdot \nabla u_h(Y,t))_{\ell_2(\Omega)}.
    \end{aligned}
\end{equation}
Now we use an identity:
$(u_h(Y,t), \partial_t\, u_h(Y,t))_{\ell_2(\Omega)} = \frac{1}{2} \partial_t (u_h(Y,t), u_h(Y,t))_{\ell_2(\Omega)} = \frac{1}{2} \partial_t \| u_h(Y,t) \|^2_{\ell_2}$
on the left-hand-side of \eqref{eq:stabilityKansa:secondStep} and the inequality \eqref{eq:stabilityestimate:discreteToContinuousInnerProduct}  
on the right-hand side of \eqref{eq:stabilityKansa:secondStep} to obtain:
\begin{equation}
    \label{eq:stabilityKansa:thirdStep}
    \frac{1}{2} \partial_t \| u_h \|^2_{\ell_2(\Omega)} \leq e_I(u_h, \bm \beta \cdot \nabla u_h) - (u_h, \bm \beta \cdot \nabla u_h)_{L_2(\Omega)},
\end{equation}
where $e_I$ is the integration error and where we also dropped the $(Y,t)$ notation for simplicity. 
The inner product on the right-hand-side of \eqref{eq:stabilityKansa:thirdStep} is skew-symmetric, which, since the advection field $\bm \beta$ is divergence-free. It implies that $(u_h, \bm \beta \cdot \nabla u_h)_{L_2(\Omega)} = \frac{1}{2} \int_{\partial\Omega} u_h^2\, \bm \beta \cdot \bm n\, ds$, where $n$ is 
the outward normal
of $\Omega$. Using that in \eqref{eq:stabilityKansa:thirdStep} and multiplying the resulting inequality by $2$ we obtain:
\begin{equation}
    \label{eq:stabilityKansa:fourthStep}
    \partial_t \| u_h \|^2_{\ell_2(\Omega)} \leq 2\, e_I(u_h, \bm \beta \cdot \nabla u_h) - 2\int_{\partial\Omega} u_h^2\, \bm \beta \cdot \bm n\, ds.
\end{equation}
We split the boundary $\partial\Omega$ into an inflow and an outflow part such that $\Omega = \Gamma_{\text{inflow}} \cup \Gamma_{\text{outflow}}$.
The two parts are defined as:
\begin{equation}
    \label{eq:stabilityKansa:inflowoutflow}
    \Gamma_{\text{inflow}}:=\{y\in\partial\Omega\, |\, \bm \beta \cdot \bm n < 0\}, \quad \Gamma_{\text{outflow}}:=\{y\in\partial\Omega\, |\, \bm \beta \cdot \bm n \geq 0\}.
\end{equation}
The boundary integral in \eqref{eq:stabilityKansa:fourthStep} is split into a sum:
$$- 2\int_{\partial\Omega} u_h^2\, \bm \beta \cdot \bm n\, ds = - 2\int_{\Gamma_{\text{inflow}}} u_h^2\, \bm \beta \cdot \bm n\, ds - 2\int_{\Gamma_{\text{outflow}}} u_h^2\, \bm \beta \cdot \bm n\, ds,$$
and then by using the sign of $\bm \beta \cdot \bm n$ from \eqref{eq:stabilityKansa:inflowoutflow} we bound
$- 2\int_{\Gamma_{\text{inflow}}} u_h^2\, \bm \beta \cdot \bm n\, ds \leq 2\int_{\Gamma_{\text{inflow}}} u_h^2\, |\bm \beta \cdot \bm n|\, ds$
and $- 2\int_{\Gamma_{\text{outflow}}} u_h^2\, \bm \beta \cdot \bm n\, ds \leq 0$. Using these two bounds in \eqref{eq:stabilityKansa:fourthStep} 
and \eqref{eq:stabilityestimate:discreteToContinuousInnerProduct} to define the $e_I$ term, we obtain
the final semi-discrete stability estimate:
\begin{equation}
    \begin{aligned}
        \label{eq:stabilityKansa:finalEstimate}
        \partial_t \| u_h \|^2_{\ell_2(\Omega)} &\leq 2\, C_{\int}\, h^{-1} h_y  \max_i \|u_h (\bm \beta\cdot \nabla u_h)\|_{L_\infty(K_i)} + 2\int_{\Gamma_{\text{inflow}}} u_h^2\, |\bm \beta \cdot \bm n|\, ds \\
        &\leq C h^{-2} h_y  \max_i \|\bm \beta\|_{L_\infty(K_i)} \max_i  \|u_h\|_{L_\infty(K_i)} + 2\int_{\Gamma_{\text{inflow}}} g^2\, |\bm \beta \cdot \bm n|\, ds, \\
    \end{aligned}
\end{equation}
where we also used an inverse inequality $\|\nabla u_h\|_{L_\infty(K_i)} \leq C_I h^{-1}\|u_h\|_{L_\infty(K_i)}$ \red{from \cite[Lemma 4.5.3]{BrennerBook}}, and defined $C=2 C_{\int} C_I$.
Thus, Kansa's RBF method has a solution of which the discrete norm is bounded in time by the inflow boundary term, and an integration error 
term of which the magnitude is possible to control by decreasing $h_y$ when $h$ is fixed.
\subsection{A semi-discrete stability estimate for RBF-PUM}
\label{sec:stabilityproperties:rbfpum}
The RBF-PUM solution $u_h$ defined in \eqref{eq:methods:rbfpum_globalsolution} and \eqref{eq:methods:rbfpum_globalsolution2} is $\mathcal C^1([0,T]) \times \mathcal C^2(\Omega)$,
since the subproblems defined on the patches are in $\mathcal C^1([0,T]) \times \mathcal C^2(\Omega_j)$, $j=1,..,N_p$, which are then blended together using the partition of unity weight functions
$w_j \in \mathcal C^2(\Omega)$.
The semi-discrete stability estimate is then precisely the same as for Kansa's RBF method stated in \eqref{eq:stabilityKansa:finalEstimate}, since the only
assumption to derive it, was a sufficient regularity of the global solution and a sufficient amount of oversampling used to discretize the PDE.

\subsection{A semi-discrete stability estimate for the RBF-FD method}
\label{sec:stabilityproperties:rbffd}
In the RBF-FD case we are dealing with a function space spanned using a set of discontinuous global cardinal functions. 
For that reason we cannot globally differentiate the solution over $\Omega$ and
can therefore not inherit the stability result from the derivation for Kansa's RBF method in Section \ref{sec:stabilityproperties:Kansa}.
However, we instead employ piecewise differentiation and make an estimate using similar techniques as in Kansa's method. 
We start to make an estimate as in \eqref{eq:stabilityKansa:secondStep}, but restrict the inner products to $K_i$ (defined in \eqref{eq:voronoi_region_Ki}):
\begin{equation}
    \label{eq:stabilityrbffd:secondStep_tmp1}
    \begin{aligned}
        (u_h(Y,t), \partial_t\, u_h(Y,t))_{\ell_2(K_i)} & = - (u_h(Y,t), \bm  \beta \cdot \nabla u_h(Y,t))_{\ell_2(K_i)},\quad i=1,..,N.
    \end{aligned}
\end{equation}
We first sum the equation above over all $K_i$ and rewrite the left-hand-side term in terms of $\Omega$ by using that $\cup_{i=1}^N K_i = \Omega$. 
We have:
\begin{equation}
    \label{eq:stabilityrbffd:secondStep}
    \begin{aligned}
        (u_h, \partial_t\, u_h)_{\ell_2(\Omega)} & = \sum_{i=1}^N - (u_h, \bm \beta \cdot \nabla u_h)_{\ell_2(K_i)},
    \end{aligned}
\end{equation}
where we also dropped the $(Y,t)$ notation.
Just as in \eqref{eq:stabilityKansa:thirdStep} we rewrite \eqref{eq:stabilityrbffd:secondStep} so that we obtain the time derivative
of a norm on the left, and use the inequality \eqref{eq:stabilityestimate:discreteToContinuousInnerProduct} 
for the integration error of piecewise continuous functions on the right, and obtain:
\begin{equation}
    \label{eq:stabilityrbffd:thirdStep}
    \partial_t \| u_h \|^2_{\ell_2(\Omega)} \leq 2\, e_I(u_h, \bm \beta \cdot \nabla u_h) - 2\sum_{i=1}^N  (u_h, \bm \beta \cdot \nabla u_h)_{L_2(K_i)}.
\end{equation}
Analogously to \eqref{eq:stabilityKansa:fourthStep}, we now integrate each inner product by parts on the right side of \eqref{eq:stabilityrbffd:thirdStep}, 
where we also use that $\nabla \cdot \bm \beta = 0$, and obtain:
\begin{equation}
    \label{eq:stabilityrbffd:fourthStep}
    \partial_t \| u_h \|^2_{\ell_2(\Omega)} \leq 2\, e_I(u_h, \bm \beta \cdot \nabla u_h) - 2\sum_{i=1}^N \int_{\partial K_i} u_h^2\, \bm \beta \cdot \bm n_{K_i}\, ds,
\end{equation}
where $\partial K_i$ is the boundary of $K_i$ as defined in \eqref{eq:voronoi_region_Ki} and $n_{K_i}$ is the outward normal of $K_i$. 
Now we take care of the Voronoi boundary terms on the right of \eqref{eq:stabilityrbffd:fourthStep}
by noting that each $\partial K_i$ is made of connected Voronoi edges. Furthermore, some $\partial K_i$ are contained inside $\Omega$ and some on $\partial \Omega$.
This implies that two neighboring interior Voronoi cells always share an edge,
where each cell has an opposite normal direction. Let us first define the set of all interior
Voronoi edges as $\mathcal{E}^I = \mathcal{E}^I_+ \cup \mathcal{E}^I_-$, where $\mathcal{E}^I_+$ are the interior edges with
a positive normal direction $\bm n^+$ and $\mathcal{E}^I_-$ the interior edges with a negative normal direction $\bm n^-$. The two normal directions 
are defined analogously to \eqref{eq:stabilityKansa:inflowoutflow} as:
\begin{equation}
    \label{eq:stabilityrbffd:normalVoronoi}
    \bm n^-_{K_i} = \{\bm n|_{\partial K_i}: \bm \beta \cdot \bm n < 0 \}, \quad \bm n^+_{K_i} = \{\bm n|_{\partial K_i}: \bm \beta \cdot \bm n \geq 0 \}.
\end{equation}
Then the sum on the right hand side of \eqref{eq:stabilityrbffd:fourthStep} is rewritten in terms of $\mathcal{E}^I$
and $\partial\Omega$:
\begin{equation}
    \begin{aligned}
    \partial_t \|u_h\|^2_{\ell_2(\Omega)} &\leq 2\, e_I(u_h, \bm \beta \cdot \nabla u_h) -  2\Big[ \int_{\partial\Omega} u_h^2\, \bm \beta \cdot \bm n\, ds\, +   \\
        & \quad \quad \qquad \quad \quad \qquad + \sum_{\mathcal{E}_i \in \mathcal{E}^I} \int_{\mathcal{E}_i} \left((u_h^+)^2\, \bm \beta \cdot \bm n^+_{\mathcal{E}_i} + (u_h^-)^2\, \bm \beta \cdot \bm n^-_{\mathcal{E}_i} \right)\, ds \Big ].
    \end{aligned}
\end{equation}
    The positive and the negative normals along the same edge $\mathcal{E}_i$ are related by $\bm n^-_{\mathcal{E}_i} = - \bm n^+_{\mathcal{E}_i}$. 
    Further, we split the integral over $\partial\Omega$ in the inflow part $\Gamma_{\text{inflow}}$ and the 
    outflow part $\Gamma_{\text{outflow}}$, both defined in \eqref{eq:stabilityKansa:inflowoutflow} such that
    $\int_{\partial\Omega} u_h^2\, \bm \beta \cdot \bm n\, ds = \int_{\Gamma_{\text{inflow}}} u_h^2\, \bm \beta \cdot \bm n\, ds + \int_{\Gamma_{\text{outflow}}} u_h^2\, \bm \beta \cdot \bm n\, ds$.
    We obtain:
    \begin{equation}
        \label{eq:stabilityrbffd:sixthstep}
        \begin{aligned}
        \partial_t \|u_h\|^2_{\ell_2(\Omega)} \leq& 2\, e_I(u_h, \bm \beta \cdot \nabla u_h) - 2\Big[ \int_{\Gamma_{\text{inflow}}} u_h^2\, \bm \beta \cdot \bm n\, ds + \int_{\Gamma_{\text{outflow}}} u_h^2\, \bm \beta \cdot \bm n\, ds\, + \\
        & \quad \quad \qquad \quad \qquad + \sum_{\mathcal{E}_i \in \mathcal{E}^I} \int_{\mathcal{E}_i} \left[(u_h^+)^2 - (u_h^-)^2\right]\, \bm \beta \cdot \bm n^+_{\mathcal{E}_i}\, ds \Big ].
        \end{aligned}
    \end{equation}
    In the same way as in Section \ref{sec:stabilityproperties:Kansa}, we now use the sign of $\bm \beta \cdot \bm n$ from \eqref{eq:stabilityKansa:inflowoutflow} 
    to further bound \eqref{eq:stabilityrbffd:sixthstep} by using $- \int_{\Gamma_{\text{inflow}}} u_h^2\, \bm \beta \cdot \bm n\, ds \leq \int_{\Gamma_{\text{inflow}}} u_h^2\, |\bm \beta \cdot \bm n|\, ds$ 
    and $- \int_{\Gamma_{\text{outflow}}} u_h^2\, \bm \beta \cdot \bm n\, ds \leq 0$, and we also use \eqref{eq:stabilityestimate:discreteToContinuousInnerProduct} to define the $e_I$ term, 
    and arrive to:
    \begin{equation}
        \label{eq:stabilityrbffd:finalEstimate}
        \begin{aligned}
        \partial_t \|u_h\|^2_{\ell_2(\Omega)} &\leq 2\, C_{\int}\, h^{-1} h_y  \max_i \|u_h (\bm \beta\cdot \nabla u_h)\|_{L_\infty(K_i)} + \\
          & \quad + 2\Big[ \int_{\Gamma_{\text{inflow}}} u_h^2\, |\bm \beta \cdot \bm n|\, ds\,  - \sum_{\mathcal{E}_i \in \mathcal{E}^I} \int_{\mathcal{E}_i} \left[(u_h^+)^2 - (u_h^-)^2\right]\, \bm \beta \cdot \bm n^+_{\mathcal{E}_i}\, ds \Big ] \\
          & \leq 2\, C h^{-2} h_y  \max_i \|\bm \beta\|_{L_\infty(K_i)} \max_i \|u_h\|_{L_\infty(K_i)} + \\
          & \quad + 2\Big[ \int_{\Gamma_{\text{inflow}}} g^2\, |\bm \beta \cdot \bm n|\, ds\,  - \sum_{\mathcal{E}_i \in \mathcal{E}^I} \int_{\mathcal{E}_i} \left[(u_h^+)^2 - (u_h^-)^2\right]\, \bm \beta \cdot \bm n^+_{\mathcal{E}_i}\, ds \Big ], \\
        \end{aligned}
    \end{equation}
    where, we used an inverse inequality \red{for general locally Sobolev-type finite dimensional spaces} \red{from \cite[Lemma 4.5.3]{BrennerBook}} in the second step in a similar way as in \eqref{eq:stabilityKansa:finalEstimate}, 
    and where we consequently have that $C=2 C_{\int} C_I$. 
    This is the final semi-discrete stability estimate for the RBF-FD method with oversampling. From \eqref{eq:stabilityrbffd:finalEstimate} 
    we observe that the solution norm is bounded by the inflow boundary term, by the jump term and by the numerical integration term. The latter two terms 
    do not originate from the physical properties of the PDE problem, however, the numerical integration error estimate term is positive and 
    can be made arbitrarily small for a fixed $h$, by decreasing $h_y$ (increasing the oversampling). 
    The jump term is more problematic because it has an arbitrary sign, and there is no obvious control over it to make it smaller. 
    \red{Our conclusion is that the jump term and the integration error term require special attention when stabilizing the RBF-FD methods in time. As is later confirmed 
    in the numerical experiments section, it is those spurious terms that make the RBF-FD method unstable in time.}

\subsection{The magnitudes of the spurious jumps in the RBF-FD cardinal functions} 
\begin{figure}[h!]
    \begin{tabular}{ccc}
        \multicolumn{3}{c}{\textbf{The RBF-FD jump magnitudes}} \\
        \hspace*{1cm}$\mathbf{p=3}$ & \hspace*{1cm}$\mathbf{p=4}$ & \hspace*{1cm}$\mathbf{p=5}$ \\
    \includegraphics[width=0.3\linewidth]{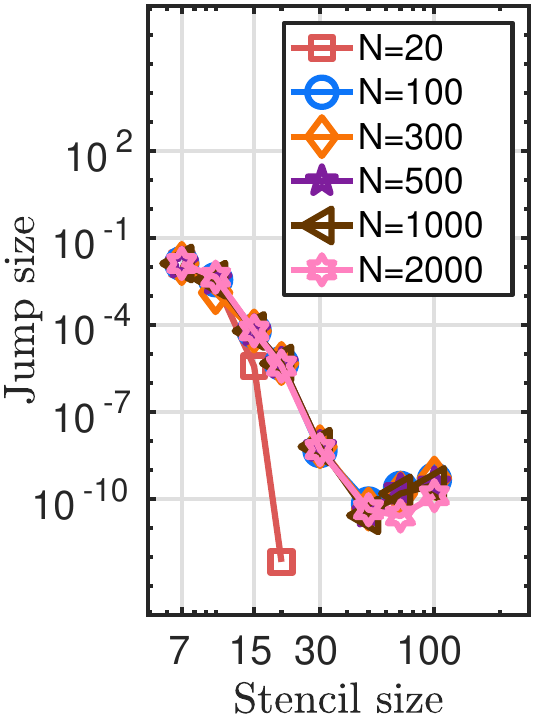} &
    \includegraphics[width=0.3\linewidth]{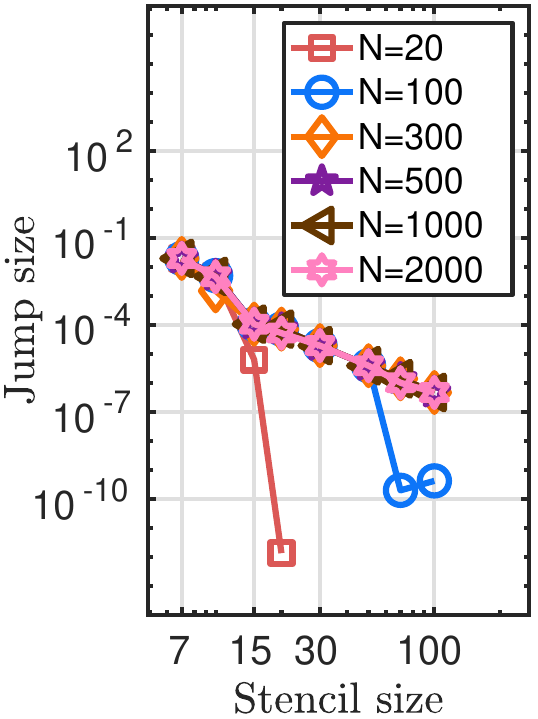} &
    \includegraphics[width=0.3\linewidth]{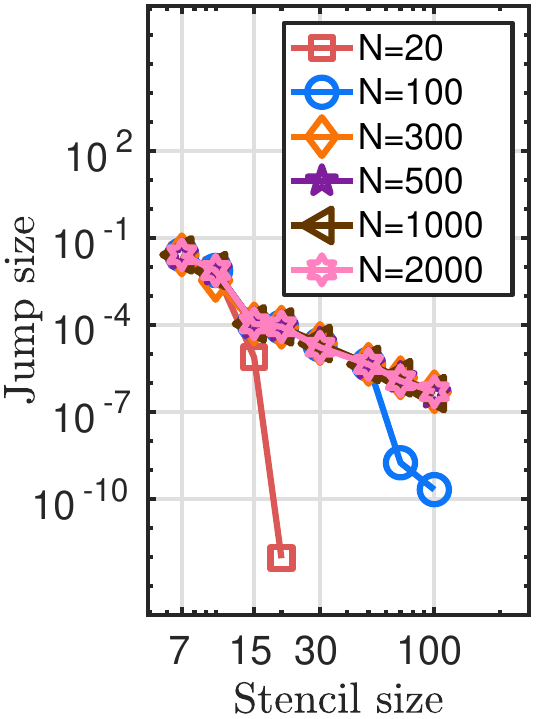} 
    \end{tabular}
    \caption{Magnitude of the largest jump (discontinuity) in a cardinal basis function $\Psi^*$, as the stencil size is increased. The 
    different curves correspond 
    to the parameter $N$, the number of nodes 
    placed on a 1D domain $\Omega_{\text{1D}}=[0,1]$. For each stencil size and $N$, the observed $\Psi^*$ is 
    chosen such that its center node is 
    closest to the point $x=0.4$.}
    \label{fig:jumps}    
\end{figure}

An a priori estimate of the magnitude of the jumps in the RBF-FD cardinal basis functions is at the present moment not available. 
Therefore we investigate these jumps numerically.
For simplicity, we focus on a 1D case, where $\Omega$ is the interval $[0,1]$. The objective is to 
understand the behavior of the discontinuities present in the RBF-FD trial space as a function of the stencil size $n$, 
the number of all nodes $N$ in the domain (proportional to $\frac{1}{h}$), and the polynomial degree $p$. 
In particular, we have a closer look at a single RBF-FD cardinal basis function $\Psi^*$, of which the center node is chosen 
closest to the point $x=0.4$. To examine the discontinuities, we measure the largest magnitude of discontinuity 
present in $\Psi^*$.
The results are given in Figure \ref{fig:jumps}, where we observe that the largest 
discontinuity decays as $n$ grows. In the cases $p=4$ and $p=5$ we observe that 
the largest discontinuity decays with the same decay rate for all tested $N$, except $N=20$ and $N=100$, where 
the discontinuity vanishes as stencil sizes become $n=20$ and $n=100$ respectively (note that $n=N$ in these two cases, and the method is in fact global). 
When $p=3$, then the largest discontinuity vanishes (up to the round-off error level) for all involved $N$. 
Another observation is that the jump size is generally the same for each $N$ when $n$ is fixed.
We draw two conclusions:
\begin{itemize}
\item For a fixed $h$, the discontinuities decrease as $n$ increases.
\item For a fixed $n$, the discontinuities remain constant as the internodal distance $h$ is decreased.
\end{itemize}
A consequence is that the jump instability inferred from \eqref{eq:stabilityrbffd:finalEstimate} can, to some extent, 
be stabilized by increasing the stencil size.

\subsection{Definition of a penalty term for the RBF-FD method that forces the spurious jumps towards zero}
\label{sec:stabilityproperties:rbffd_penalty}
\begin{figure}[h!]
    \centering
    \includegraphics[width=0.23\linewidth]{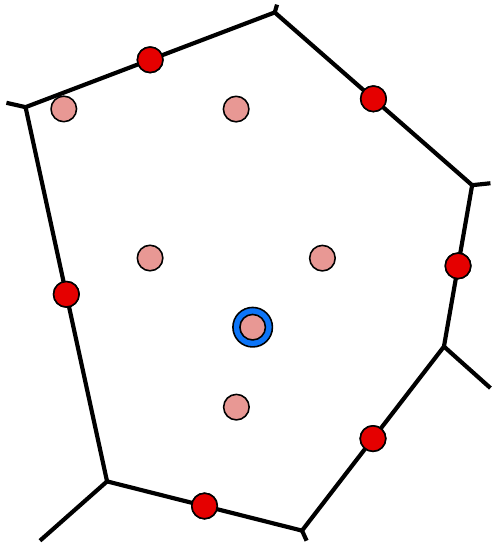}
\caption{One Voronoi cell with a set of midpoints $Y_{\mathcal{E}}$ placed over its edges (red points). 
The Voronoi cell has one center $X$ point (blue point) and 6 interior evaluation $Y$ points (pale red points).}
\label{fig:jumps:stabilization_midpoints}    
\end{figure}

A common way to stabilize the system 
\eqref{eq:methods:advection_system_matrix_projected}, where the matrices are constructed using the RBF-FD method, 
is to supply \eqref{eq:methods:advection_system_matrix_projected} with a stabilization term $\gamma \mathcal{P} u_h$ on the right-hand-side of the system. 
Here $\mathcal{P}$ is normally a
hyperviscosity operator and $\gamma$ is a scaling parameter. In this section we provide an alternative formulation of $\mathcal{P}$ based on 
the stability estimate \eqref{eq:stabilityrbffd:finalEstimate} where the jumps across the Voronoi edges were shown 
to cause a spurious growth in time.
For this reason we construct $\mathcal{P}$ such that the jump terms $\sum_{\mathcal{E}_i \in \mathcal{E}^I}\int_{\mathcal{E}_i} (u_h^+-u_h^-)^2\, ds$ are forced to $0$. 
The discrete operator $\mathcal{P}$ and the scaling $\gamma$ are given by discretizing the jump term using the oversampled matrices:
\begin{eqnarray*}
    \mathcal{P} &=& \left(E_+(Y_\mathcal{E},X) - E_-(Y_\mathcal{E},X)\right)^T \left(E_+(Y_\mathcal{E},X) - E_-(Y_\mathcal{E},X)\right),\\
    \gamma &=& -h_{\mathcal{E}},\,
\end{eqnarray*}
where $Y_\mathcal{E}$ is the set of all midpoints on the interior Voronoi edges, see Figure \ref{fig:jumps:stabilization_midpoints}.
The terms $E_+(Y_\mathcal{E},X)$ and $E_-(Y_\mathcal{E},X)$ are unscaled global matrices of evaluation weights for evaluating the solution at 
$Y_\mathcal{E}$, 
constructed using the ''left'' stencil and the ''right'' stencil respectively. 
The $h_{\mathcal{E}}>0$ scaling is the approximate mean distance between the points $Y_\mathcal{E}$, measured along the Voronoi edges.

\section{Numerical experiments}
\label{sec:experiments}
In this section we provide numerical tests that (i) show a practical relation to the derived \red{semi-discrete} stability estimates in 
Section \ref{sec:stabilityproperties:Kansa} for Kansa's method, 
in Section \ref{sec:stabilityproperties:rbfpum} for RBF-PUM 
and in Section \ref{sec:stabilityproperties:rbffd} for the RBF-FD method, (ii) examine the convergence of the approximation error under node refinement. 
\red{In all experiments we use the Runge-Kutta 4 method when discretizing the solution to the advection problem in time.}
The domain boundary $\partial \Omega$ that we use for the experiments is in polar coordinates $(r,\theta)$ prescribed by the relation 
$r(\theta) = 1 - \frac{1}{3} \sin(2\theta)^2$, where $\theta \in [0, 2\pi)$. The final shape can be observed in Figure \ref{fig:experiments:pointset}. 

The initial condition that we consider is a $\mathcal C^6(\Omega)$ compactly supported Wendland function $(1-r)^6\, (35 r^2 + 18r + 3)$, where $r$ is the Euclidean distance measure from the origin $(0,0)$.  
The initial condition is scaled such that its support radius is of size $0.4$. 

Time steps for advancing the numerical solution in time and scaling the eigenvalue spectra, 
are computed using the relation:
\begin{equation}
    \label{eq:experiments:timestep}
\Delta t_n = \text{CFL} \frac{h}{\|\bm \beta(\cdot, t_n) \|_{L_{\infty}(\Omega)}},
\end{equation}
where $\text{CFL}>0$ is the CFL (Courant-Friedrichs-Lax) number, $\|\bm \beta (\cdot, t_n) \|_{L_{\infty}(\Omega)}$ is the maximum wave-speed over the whole domain at time $t=t_n$. 
We consider a rotational field $\bm \beta = \frac{1}{2} \left [\cos(2\pi t), \sin(2\pi t) \right ]$. 

The approximation errors are measured as:
\begin{equation}
    \label{eq:experiments:errors}
    \|e\|_k = \frac{\|u_h(Y,t_{\text{final}}) - u(Y,t_{\text{final}})\|_k}{\|u(Y,t_{\text{final}})\|_k},\quad k=\{1,2,\infty\}.
\end{equation}

The point set $X$ is obtained using the DistMesh algoritm \cite{DistMesh} with $h$ as input, and is then in some cases modified with a random perturbation 
of the magnitude $0.65\, h$ in order to also make observations on 
a scattered point set. A comparison between the two instances of point sets when $h=0.09$, is given in Figure \ref{fig:experiments:pointset}. 
In the table below we give a relation between all considered $h$ and the number of unknowns $N$.
\begin{center}
    \vspace{0.2cm}
\begin{tabular}{|c|c|c|c|c|c|c|c|}
    \hline     
    $h$ & 0.08 & 0.07 & 0.06 & 0.05 & 0.04 & 0.03 & 0.02 \\ \hline
    $1/h$ & 12.5  & 14.3 &  16.7 &  20 &  25 &  33.3 &  50 \\ \hline
    $N$ & 400 &  520 & 713 & 1030 & 1605 & 2850 & 6420 \\
    \hline 
\end{tabular}
\,     \vspace{0.2cm}
\end{center}
The point set $Y$ is in the cases when we examine the solutions to the PDE problem \eqref{eq:intro:conservationlaw}, 
obtained using the DistMesh algorithm.

\subsection{Convergence of the integration error induced by oversampling}
\label{sec:experiments:integration_error}
Here we numerically verify the asymptotic behavior of the integration error estimate \eqref{appendix:integral:Omega} derived in Appendix \ref{sec:appendix:integration}. 
We observe the convergence 
$\frac{|\Omega|}{M} \sum_{k=1}^M f(y_k) \to \int_\Omega f(y)\, d\Omega,$ where $|\Omega| = \int_\Omega 1\, d\Omega$, 
and the evaluation points $Y = \{y_k\}_{k=1}^M$ are chosen as: Halton points, Cartesian points or DistMesh points.
The domain boundary that we consider 
is defined in the scope of Section \ref{sec:experiments}, and is visualized in Figure \ref{fig:experiments:pointset}. 


Let $r_f = r_f(y_1, y_2) = (y_1^2 + y_2^2)^{1/2}$. 
We choose three functions that we integrate: 
a Gaussian function $f_1(y_1, y_2) = e^{-3 r_f}$ (infinitely smooth), 
a cubic polyhanormic spline $f_2(y_1, y_2) = r_f^3 $ (twice continuously differentiable),
a discontinuous function defined by $f_3(y_1,y_2)=0.2 + \sin(4\pi y_1 y_2)$ when $r_f\leq 0.5$, 
$f_3(y_1,y_2)=y_1^3 y_2$ when $0.5 < r_f \leq 0.7$ and $f_3(y_1,y_2) = 0.4 + \cos(4\pi y_1 y_2)$ when $r_f > 0.7$.

\begin{figure}
\begin{tabular}{ccc}
        \multicolumn{3}{c}{\hspace{0.2cm} \textbf{Integration error convergence when integrating 3 functions}} \\
\hspace{1cm} \textbf{Gaussian} & \hspace{0.7cm} \textbf{Cubic PHS} & \hspace{0.7cm} \textbf{Discontinuous} \\
\includegraphics[width=0.3\linewidth]{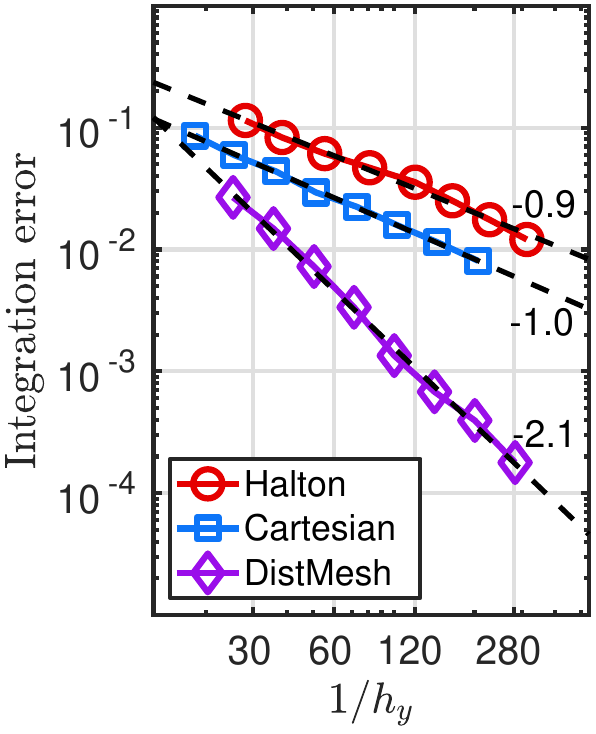} &
\includegraphics[width=0.3\linewidth]{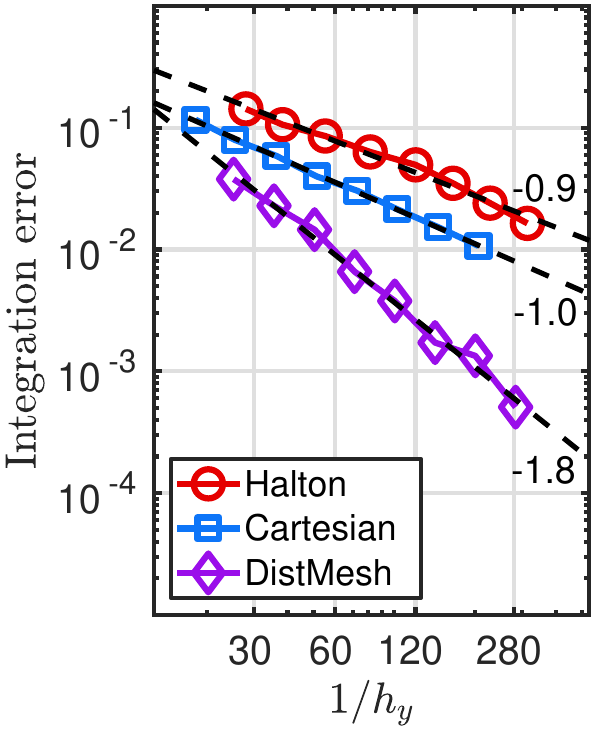} &
\includegraphics[width=0.3\linewidth]{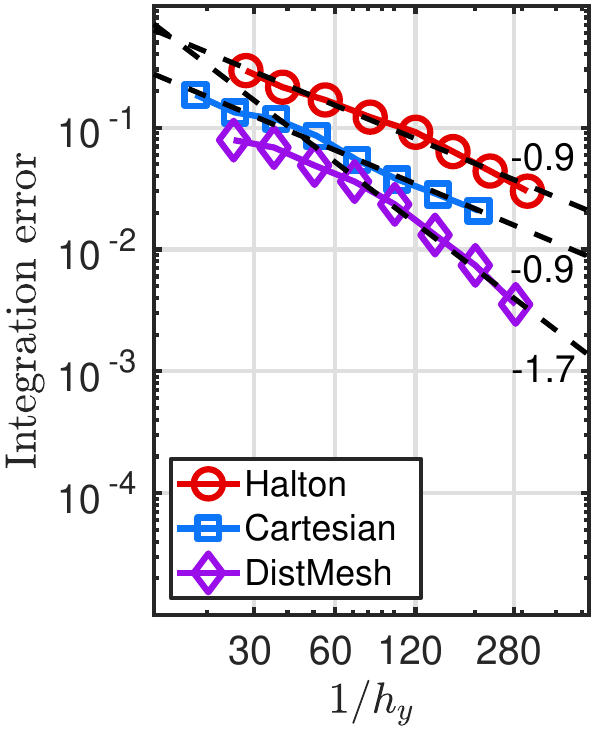}     
\end{tabular}    
\caption{Integration error when approximating an integral by means of oversampling, as a function of the inverse internodal distance $1/h_y$ in the evaluation point set $Y$. 
The three plots are from left to right displaying the convergence when integrating: a Gaussian function, the cubic polyharmonic (PHS) function and a discontinuous function 
disjointly defined by trigonometric and polynomial functions. The lines in each plot correspond to using three different evaluation point sets $Y$.}
\label{fig:experiments:integration:domain_1}
\end{figure}
The result is given in Figure \ref{fig:experiments:integration:domain_1}, where we observe that the convergence trend is $1$ for all considered integrands, 
when using Halton and Cartesian point sets for evaluating the integrals. We also observe that the Cartesian point set in the considered case gives 
a smaller integration error. The convergence trend when the DistMesh point set is used tends towards $2$. 
The most significant difference between the point sets is that the internodal distance near the domain boundary is following a smooth distribution 
in the DistMesh point set case, and we speculate that this is the reason for observing a larger convergence trend. In any case, the integration error 
due to the oversampling decays asymptotically with at least an order $1$ in 2D, as derived in \eqref{appendix:integral:Omega}.

We note that we observed very similar convergence trends of the integration error when integrating the three considered functions over a different 
choice of domain boundary, in polar coordinates 
given by $r(\theta) = 1+1/10(\sin(7t)+\sin(t))$, $\theta \in [0, 2\pi]$. A visualization of that domain boundary is available in \cite[Fig.~3]{ToLaHe21}. 
\subsection{Eigenvalue spectra for Kansa's RBF method, RBF-PUM and the RBF-FD method}
\label{sec:experiments:eigenvalues}
When computing the eigenvalue spectra we impose a $0$ inflow boundary condition by 
eliminating the corresponding columns and rows of $(E_h^T E_h)^{-1}\, (E_h^T\, D_h)$.
For fixed parameters $h$ and $p$ 
we test how an increase in the oversampling parameter $q$ influences the stability 
through the eigenvalue spectra. Note that as $h_y \to 0$, then $q \to \infty$, and the integration error in the stability estimate 
\eqref{eq:stabilityKansa:finalEstimate} goes to $0$. Note also that when $q=1$, we have that $h_y=h$, which is the limit case of the oversampling 
that yields the  collocated discretization.

\begin{figure}
    \centering
\begin{tabular}{ccccccccc}
    \multicolumn{9}{c}{\textbf{Kansa's method: eigenvalue spectra, DistMesh nodes} \vspace{0.2cm}} \\
\includegraphics[height=0.35\linewidth]{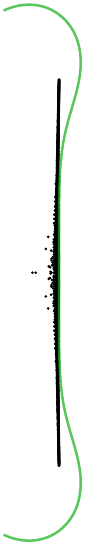} &
\includegraphics[height=0.35\linewidth]{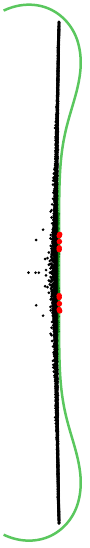} &
\includegraphics[height=0.35\linewidth]{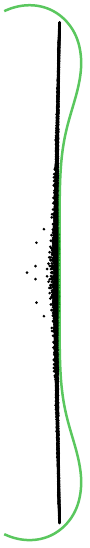} &
\includegraphics[height=0.35\linewidth]{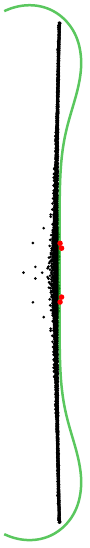} &
\includegraphics[height=0.35\linewidth]{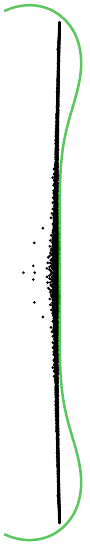} &
\includegraphics[height=0.35\linewidth]{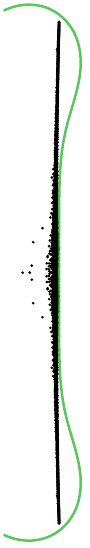} &
\includegraphics[height=0.35\linewidth]{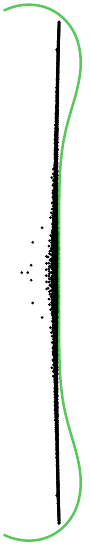} &
\includegraphics[height=0.35\linewidth]{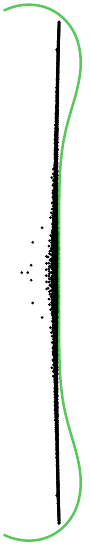} &
\includegraphics[height=0.35\linewidth]{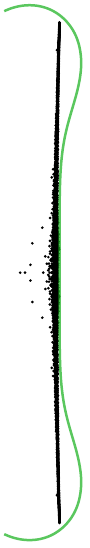} \\
$\mathbf{q=1}$ & $\mathbf{q=2}$ & $\mathbf{q=3}$ & $\mathbf{q=4}$ & $\mathbf{q=5}$ & $\mathbf{q=6}$ & $\mathbf{q=7}$ & $\mathbf{q=8}$ & $\mathbf{q=9}$  
\end{tabular}
\caption{Eigenvalue spectra (black dots) of the advection operator with $0$ inflow boundary conditions, when Kansa's method is used, as the oversampling parameter $q$ is increased (left to right). 
The green line is the boundary of the classical explicit Runge-Kutta 4 stability region. The red dots are the eigenvalues that are not contained inside the stability region.}
\label{fig:experiments:eigenvalues:Kansa_unperturbed}
\end{figure}
\begin{figure}
    \centering
\begin{tabular}{ccccccccc}
    \multicolumn{9}{c}{\textbf{Kansa's method: eigenvalue spectra, perturbed DistMesh nodes} \vspace{0.2cm}} \\
\hspace{-0.3cm}\includegraphics[width=0.2\linewidth]{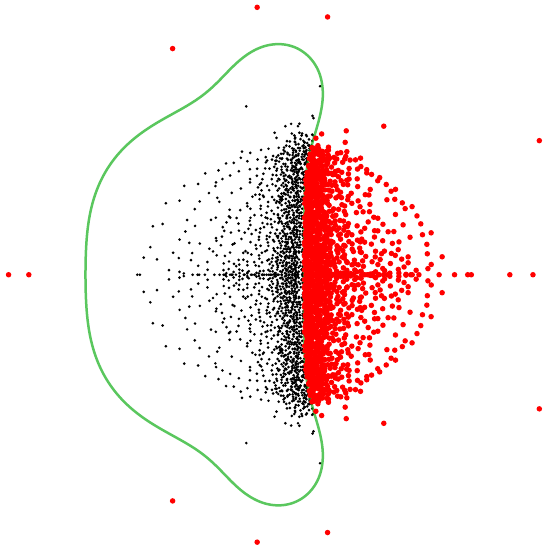} & 
\hspace{-0.1cm}\includegraphics[height=0.35\linewidth]{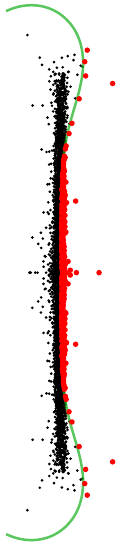} &
\hspace{-0.1cm}\includegraphics[height=0.35\linewidth]{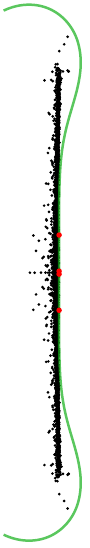} &
\hspace{-0.1cm}\includegraphics[height=0.35\linewidth]{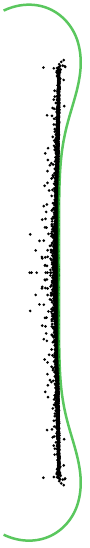} &
\hspace{-0.1cm}\includegraphics[height=0.35\linewidth]{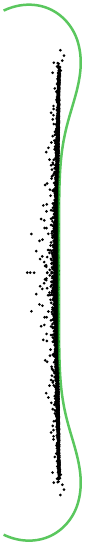} &
\hspace{-0.1cm}\includegraphics[height=0.35\linewidth]{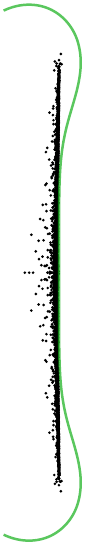} &
\hspace{-0.1cm}\includegraphics[height=0.35\linewidth]{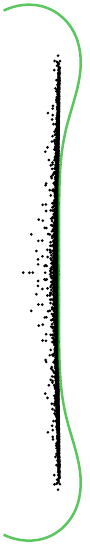} &
\hspace{-0.1cm}\includegraphics[height=0.35\linewidth]{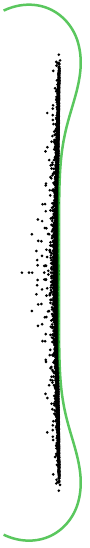} &
\hspace{-0.1cm}\includegraphics[height=0.35\linewidth]{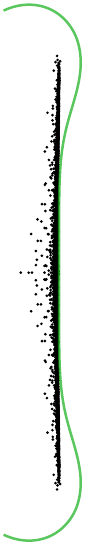}\\

\hspace{-0.3cm} $\mathbf{q=1}$ & \hspace{-0.5cm} $\mathbf{q=2}$ & \hspace{-0.2cm} $\mathbf{q=3}$ & \hspace{-0.2cm} $\mathbf{q=4}$ & \hspace{-0.2cm} $\mathbf{q=5}$ & \hspace{-0.2cm} $\mathbf{q=6}$ & \hspace{-0.2cm} $\mathbf{q=7}$ & \hspace{-0.2cm} $\mathbf{q=8}$ & \hspace{-0.2cm} $\mathbf{q=9}$ 
\end{tabular}
\caption{Eigenvalue spectra (black dots) of the advection operator with $0$ inflow boundary conditions, 
when Kansa's method is used on randomly perturbed DistMesh $X$ points ($h=0.03$), 
as the oversampling parameter $q$ is gradually increased. 
The green line is the full stability region ($q=1$) or its boundary ($q>1$), of the classical explicit Runge-Kutta 4 method. 
The red dots are the eigenvalues that are not contained inside the stability region.}
\label{fig:experiments:eigenvalues:Kansa_perturbed}
\end{figure}

The result for Kansa's method with the DistMesh $X$ points when $h=0.03$ and $p=4$, is drawn in Figure \ref{fig:experiments:eigenvalues:Kansa_unperturbed}. 
We observe that the spectra are alternating between being stable and having spurious eigenvalues, until $q=4$, and then 
become stable for $q\geq 5$. A similar result is also observed when the $X$ points are randomly perturbed in Figure \ref{fig:experiments:eigenvalues:Kansa_perturbed}. 
Both results are aligned with the stability result from Section \ref{sec:stabilityproperties:Kansa}: for an oversampling that is large enough, 
Kansa's method is stable when discretizing the linear advection problem.

\begin{figure}
    \centering
\begin{tabular}{ccccccccc}
    \multicolumn{9}{c}{\textbf{RBF-PUM: eigenvalue spectra, perturbed DistMesh nodes} \vspace{0.2cm}} \\
\hspace{-0.2cm}\includegraphics[width=0.23\linewidth]{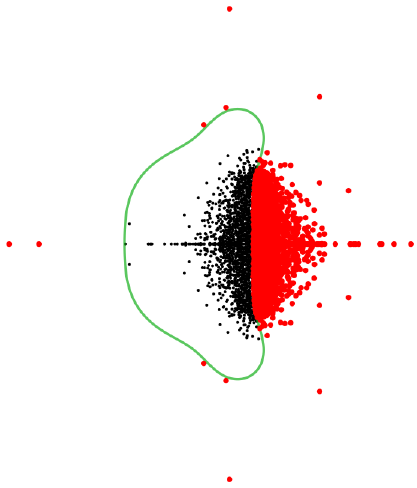} &
\hspace{-0.4cm}\includegraphics[height=0.35\linewidth]{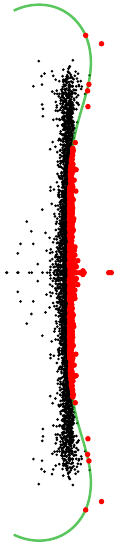} &
\hspace{-0.2cm}\includegraphics[height=0.35\linewidth]{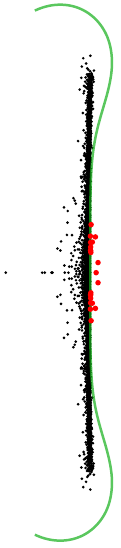} &
\hspace{-0.2cm}\includegraphics[height=0.35\linewidth]{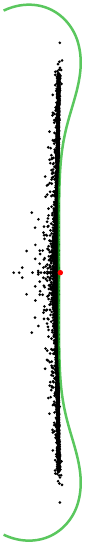} &
\hspace{-0.2cm}\includegraphics[height=0.35\linewidth]{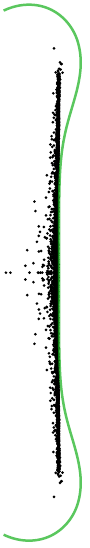} &
\hspace{-0.2cm}\includegraphics[height=0.35\linewidth]{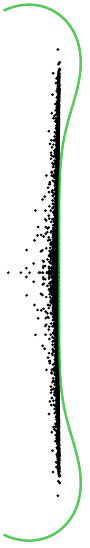} &
\hspace{-0.2cm}\includegraphics[height=0.35\linewidth]{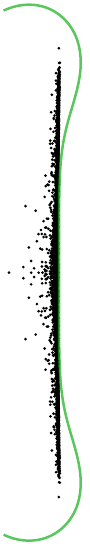} &
\hspace{-0.2cm}\includegraphics[height=0.35\linewidth]{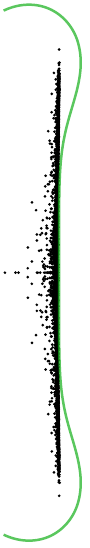} &
\hspace{-0.1cm}\includegraphics[height=0.35\linewidth]{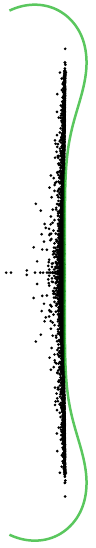}\\

$\mathbf{q=1}$ & \hspace{-0.5cm}$\mathbf{q=2}$ & \hspace{-0.1cm}$\mathbf{q=3}$ & \hspace{-0.1cm}$\mathbf{q=4}$ & \hspace{-0.1cm}$\mathbf{q=5}$ & \hspace{-0.1cm}$\mathbf{q=6}$ & \hspace{-0.1cm}$\mathbf{q=7}$ & \hspace{-0.1cm}$\mathbf{q=8}$ & \hspace{-0.1cm}$\mathbf{q=9}$
\end{tabular}
\caption{Eigenvalue spectra (black dots) of the advection operator with $0$ inflow boundary conditions, 
when RBF-PUM is used on randomly perturbed DistMesh $X$ points ($h=0.03$), 
as the oversampling parameter $q$ is gradually increased. 
The green line is the full stability region ($q=1$) or its boundary ($q>1$), of the classical explicit Runge-Kutta 4 method. 
The red dots are the eigenvalues that are not contained inside the stability region.}
\label{fig:experiments:eigenvalues:rbfpum_perturbed}
\end{figure}

The result for RBF-PUM when $h=0.03$ and $p=4$ is given in Figure \ref{fig:experiments:eigenvalues:rbfpum_perturbed}. 
Here we only consider the randomly perturbed DistMesh point set $X$. We observe that the eigenvalue spectra are stable when $q \geq 5$, but otherwise 
unstable. This is aligned with the stability result from Section \ref{sec:stabilityproperties:rbfpum}.

\begin{figure}
    \centering
\begin{tabular}{lcccccccc}
    \multicolumn{9}{c}{\textbf{The unstabilized RBF-FD method: eigenvalue spectra} \vspace{0.3cm}} \\
    \rotatebox{90}{\hspace{2cm} $\mathbf{n=12}$} &
    \includegraphics[height=0.35\linewidth]{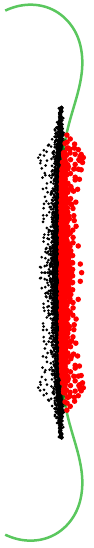} &
    \includegraphics[height=0.35\linewidth]{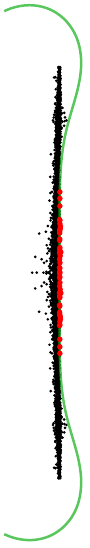} &
    \includegraphics[height=0.35\linewidth]{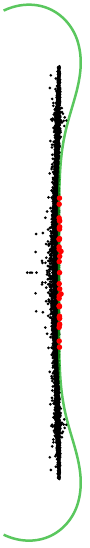} &
    \includegraphics[height=0.35\linewidth]{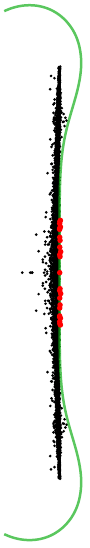} &
    \includegraphics[height=0.35\linewidth]{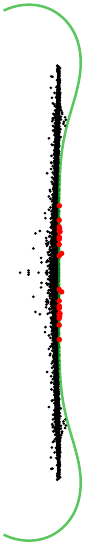} &
    \includegraphics[height=0.35\linewidth]{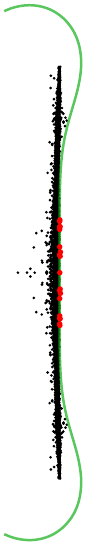} &
    \includegraphics[height=0.35\linewidth]{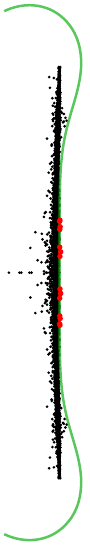} &
    \includegraphics[height=0.35\linewidth]{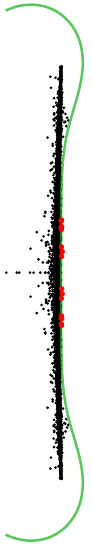}\\

    \rotatebox{90}{\hspace{2cm} $\mathbf{n=30}$} &
    \includegraphics[height=0.35\linewidth]{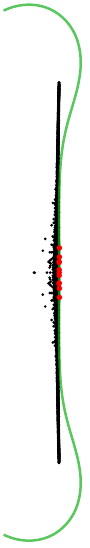} &
    \includegraphics[height=0.35\linewidth]{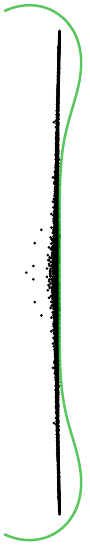} &
    \includegraphics[height=0.35\linewidth]{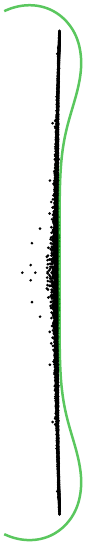} &
    \includegraphics[height=0.35\linewidth]{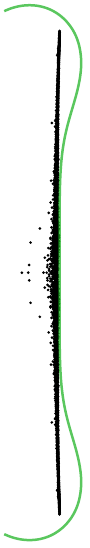} &
    \includegraphics[height=0.35\linewidth]{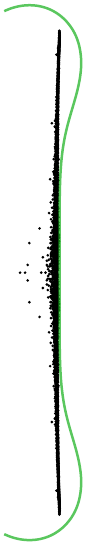} &
    \includegraphics[height=0.35\linewidth]{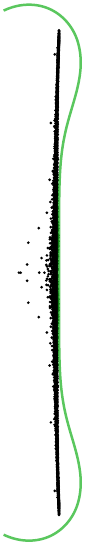} &
    \includegraphics[height=0.35\linewidth]{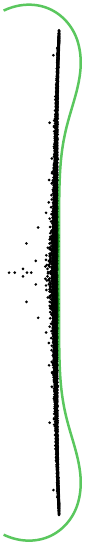} &  
    \includegraphics[height=0.35\linewidth]{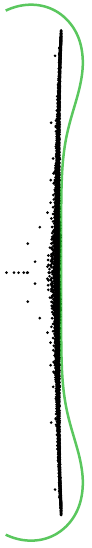}\\

& $\mathbf{q=1}$ & $\mathbf{q=3}$ & $\mathbf{q=5}$  & $\mathbf{q=7}$ & $\mathbf{q=9}$ & $\mathbf{q=10}$ & $\mathbf{q=30}$ & $\mathbf{q=60}$
\end{tabular}
\caption{Eigenvalue spectra (black dots) of the advection operator with $0$ inflow boundary conditions, when the RBF-FD method is used on DistMesh $X$ points ($h=0.03$), 
as the oversampling parameter $q$ (columns) is gradually increased. The monomial basis degree is fixed to $p=2$. 
The two rows correspond to a different stencil size $n$.
The green line is the boundary of the stability region for the classical explicit Runge-Kutta 4 method. 
The red dots are the eigenvalues that are not contained inside the stability region.}
\label{fig:experiments:eigenvalues:rbffd}
\end{figure}

The results for the unstabilized RBF-FD method with $h=0.03$ are given in Figure \ref{fig:experiments:eigenvalues:rbffd}. 
 Here we use an unperturbed DistMesh point set $X$. 
We test the eigenvalue spectra as the stencil size $n$ varies and the monomial basis degree $p=2$ is fixed.
When $n=12$ (obtained from $n = 2\binom{p+2}{2}$), the eigenvalue spectra contain spurious eigenvalues for each choice of the oversampling parameter, 
even when the oversampling is unusually large. 
This is in agreement with the semi-discrete stability estimate \eqref{eq:stabilityrbffd:finalEstimate}. 
The jump term in the estimate has a non-negligible size 
and an arbitrary sign. 
From the same figure we further observe 
that when we keep $p=2$, but increase the stencil size to $n=30$, the eigenvalue spectra become stable for $q \geq 3$. This is again in agreement with 
the estimate \eqref{eq:stabilityrbffd:finalEstimate},  when the estimate is augmented with the result from Figure 
\eqref{fig:jumps} which states that the magnitude of the jumps vanishes as $n$ is increased. Thus, large stencil sizes improve 
the stability properties of the RBF-FD method. 
We make a remark that the eigenvalue spectra were stable for larger $p$ and $n$ computed 
according to $n=2\binom{p+2}{2}$. This makes sense, since larger $p$ are associated with larger $n$, and thus smaller jumps. However, 
we can not say in general that the RBF-FD method is stable for large $p$, since the size of the jumps will inevitably depend 
on the choice of the point set $X$.

\begin{figure}
    \centering
\begin{tabular}{lccccccc}
    \multicolumn{8}{c}{\textbf{Jump-stabilized RBF-FD method: eigenvalue spectra} \vspace{0.3cm}} \\
    \rotatebox{90}{\hspace{2cm} $\mathbf{n=12}$} &
    \hspace{-0.8cm}\includegraphics[height=0.3\linewidth]{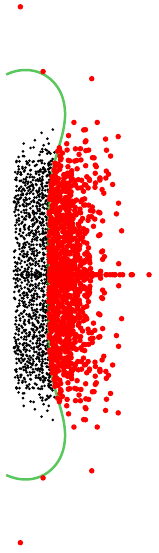} &
    \includegraphics[height=0.3\linewidth]{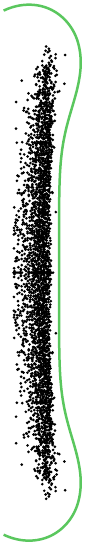} &
    \includegraphics[height=0.3\linewidth]{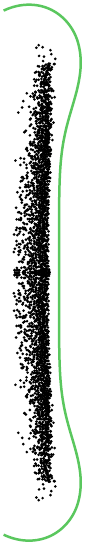} &
    \includegraphics[height=0.3\linewidth]{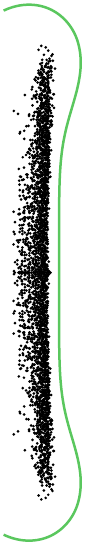} &
    \includegraphics[height=0.3\linewidth]{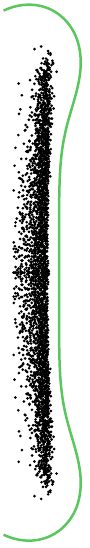} &
    \includegraphics[height=0.3\linewidth]{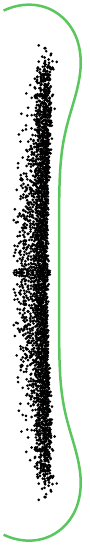} &
    \includegraphics[height=0.3\linewidth]{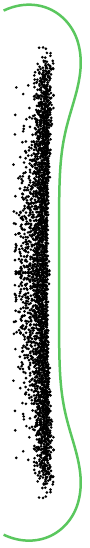} \\

    \rotatebox{90}{\hspace{2cm} $\mathbf{n=30}$} &
    \includegraphics[height=0.3\linewidth]{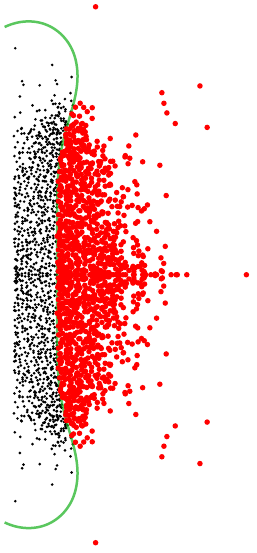} &
    \includegraphics[height=0.3\linewidth]{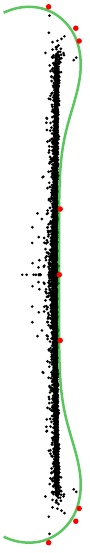} &
    \includegraphics[height=0.3\linewidth]{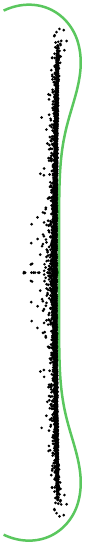} &
    \includegraphics[height=0.3\linewidth]{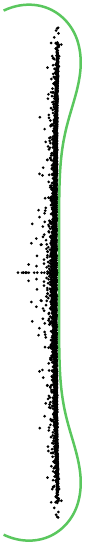} &
    \includegraphics[height=0.3\linewidth]{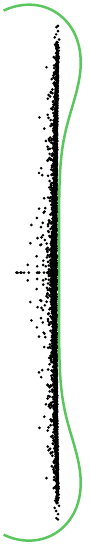} &
    \includegraphics[height=0.3\linewidth]{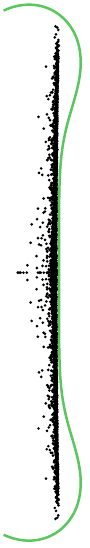} &
    \includegraphics[height=0.3\linewidth]{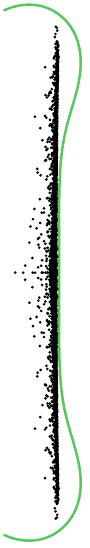} \\

& $\mathbf{q=1}$ & $\mathbf{q=3}$ & $\mathbf{q=5}$  & $\mathbf{q=7}$ & $\mathbf{q=9}$ & $\mathbf{q=10}$ & $\mathbf{q=30}$ 
\end{tabular}
\caption{Eigenvalue spectra (black dots) of the advection operator with $0$ inflow boundary conditions, when the RBF-FD method is stabilized using the jump penalty term and 
used in combination with randomly perturbed DistMesh $X$ points ($h=0.03$), 
as the oversampling parameter $q$ (columns) is gradually increased. The monomial basis degree is fixed to $p=2$. 
The two rows correspond to a different stencil size $n$.
The green line is the boundary of the stability region for the classical explicit Runge-Kutta 4 method. 
The red dots are the eigenvalues that are not contained inside the stability region.}
\label{fig:experiments:eigenvalues:rbffd_jumpstabilized}
\end{figure}

The results for the jump-stabilized RBF-FD method are collected in Figure \ref{fig:experiments:eigenvalues:rbffd_jumpstabilized}. 
We again use $p=2$ and two different stencil sizes, $n=12$ and $n=30$, but 
now we add the jump penalty term defined in Section \ref{sec:stabilityproperties:rbffd_penalty} that forces the spurious jumps in the 
RBF-FD cardinal basis functions towards $0$. We also use the randomly perturbed DistMesh point set in order to not give the stabilization term 
a favorable starting point. From the figure we observe that the eigenvalue spectrum around the imaginary axis is now stable 
also in the $n=12$ case, as opposed to the unstabilized method in Figure \ref{fig:experiments:eigenvalues:rbffd}. 
We note that we also observed stable eigenvalue spectra for higher choices of $p$ when the jump stabilization term was used and we chose $n$ 
according to the 
standard formula $n=2\binom{p+2}{2}$. Another remark is that the $\text{CFL}$ numbers have to be chosen smaller when the jumps are large 
(small $p$ and consequently small $n$), compared with the case when the RBF-FD method is not stabilized using the jump penalty term. 
Thus, an effective strategy to stabilize the RBF-FD method is to always add the jump penalty term and possibly, at the same time, 
increase the stencil size when $p$ is small ($p=1,2$).

\subsection{Convergence of the approximation error under node refinement for RBF-PUM and the RBF-FD method}
\label{sec:experiments:error_href}
We examine the convergence of the approximation error as a function of $h$, for different choices of $p$ and a fixed $q$. 
The final time of the simulation is set to $t_{\text{final}}=3$. The time step is computed according to $\text{CFL}=0.2$, which is intentionally kept low, 
so that the resulting time steps are small and that the spatial approximation errors dominate.

\begin{figure}
    \centering
\begin{tabular}{ccc}
    \multicolumn{3}{c}{\hspace{0.8cm}\textbf{Convergence under node refinement in $2$-norm} \vspace{0.1cm}} \\
    \hspace{1cm}\textbf{RBF-PUM} & \hspace{0.6cm}\textbf{RBF-FD} & \hspace{0.3cm}\textbf{RBF-FD, stabilized}\\
    \includegraphics[width=0.3\linewidth]{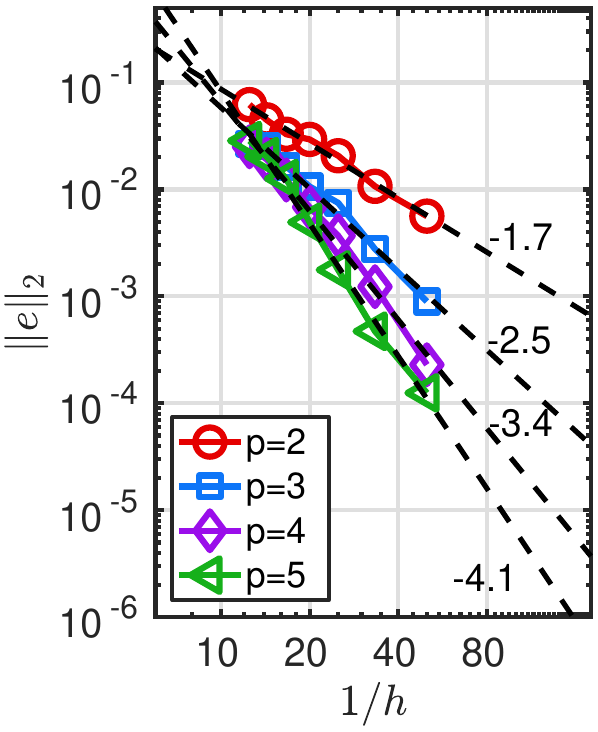} & \hspace{-0.5cm} 
    \includegraphics[width=0.3\linewidth]{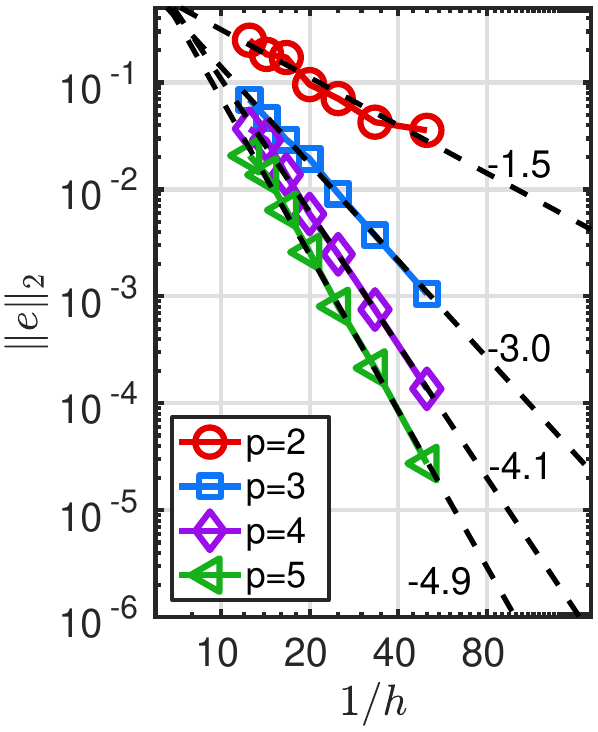} & \hspace{-0.8cm}
    \includegraphics[width=0.3\linewidth]{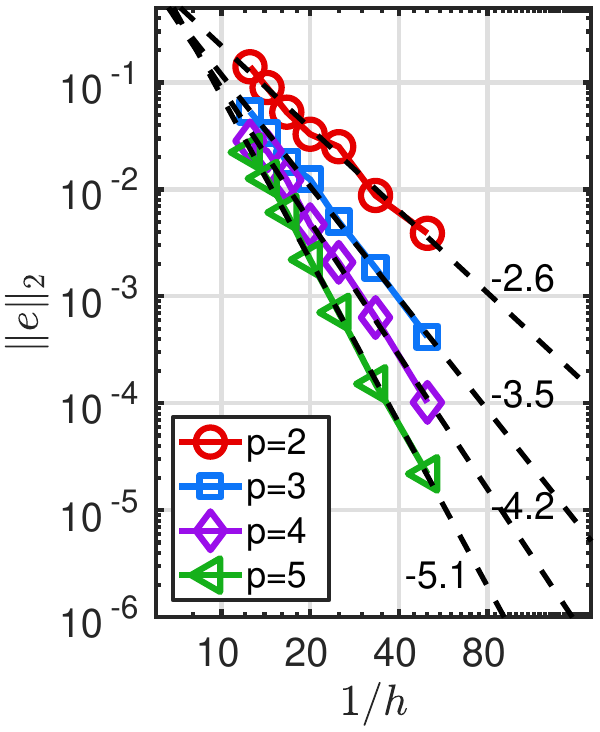} 
\end{tabular}
\caption{Convergence of the approximation error in the relative $2$-norm, as a function of the inverse internodal distance $1/h$, for three methods: RBF-PUM, unstabilized RBF-FD and jump-stabilized RBF-FD. 
Different lines in the plots correspond to different choices of the monomial basis degree $p$ used to construct the stencil-based approximations. The oversampling parameter 
is fixed to $q=9$. The time step to advance the solution in time is computed using $\text{CFL}=0.2$.}
\label{fig:experiments:error_href}
\end{figure}
The results for RBF-PUM, the unstabilized RBF-FD method and the jump-stabilized RBF-FD method, using $q=9$, are shown in Figure \ref{fig:experiments:error_href}. 
We observe that all the methods converge for all $p$ with a rate at least $p-1$. In the unstabilized RBF-FD case, 
the eigenvalue spectra for $p=2$ were shown to be unstable in Figure 
\ref{fig:experiments:eigenvalues:rbffd}, however, the approximation error was still convergent. This is because the spurious eigenvalues have not triggered 
a visible instability until $t_{\text{final}}=3$. A longer simulation that shows an instability is performed in the section that follows. 
The approximation error in the jump-stabilized RBF-FD case is at least as small 
as in the unstabilized RBF-FD case. We conclude that the proposed jump-stabilization is effective for the considered simulation case. 

Similar convergence rates were also observed in $1$-norm and $\infty$-norm, for all considered methods.

\subsection{Energy of the numerical solution as a function of time for Kansa's RBF method, RBF-PUM and the RBF-FD method}
\label{sec:experiments:conservation}
Here we examine the energy of the numerical solution with respect to time, for all three methods. 
The purpose is to examine the stability of the methods, but also to observe whether the energy is conserved. 
We measure the energy of the numerical solution relative to the energy of the initial condition as $\|u_h(Y,t)\|^2_{\ell_2}/\|u_h(Y,0)\|^2_{\ell_2}$, 
where the norm is defined in \eqref{eq:stabilityestimate:discretenorms}. 
The final time of the simulation is set to $t_{\text{final}}=20$, which makes the compactly supported initial 
condition rotate around the origin $20$ times, without interacting with the boundary of the domain, 
so that the energy of the solution is (ideally) not changing with time. The CFL number is set to $\text{CFL}=0.6$. 
We test the conservation for different choices of $h$, $p$ and $q$.

\begin{figure}
    \centering
\begin{tabular}{lccc}
    \multicolumn{4}{c}{\hspace{0.8cm}\textbf{Kansa's RBF method: energy in time} \vspace{0.1cm}} \\
    & \hspace{0.8cm} $\mathbf{q=6}$ & \hspace{0.5cm}$\mathbf{q=9}$ & \hspace{0.5cm} $\mathbf{q=12}$\\
    \rotatebox{90}{\hspace{2.2cm}$\mathbf{p=5}$}  & \includegraphics[width=0.3\linewidth]{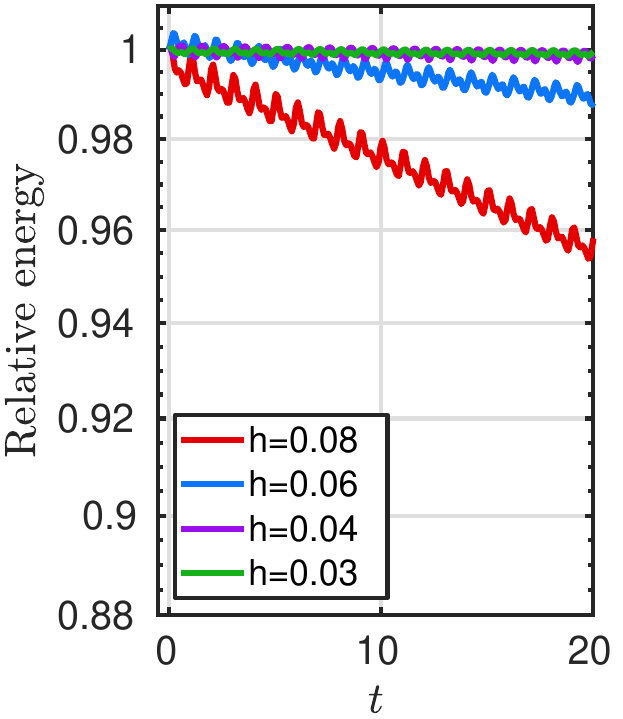} & 
    \hspace{-0.4cm}\includegraphics[width=0.3\linewidth]{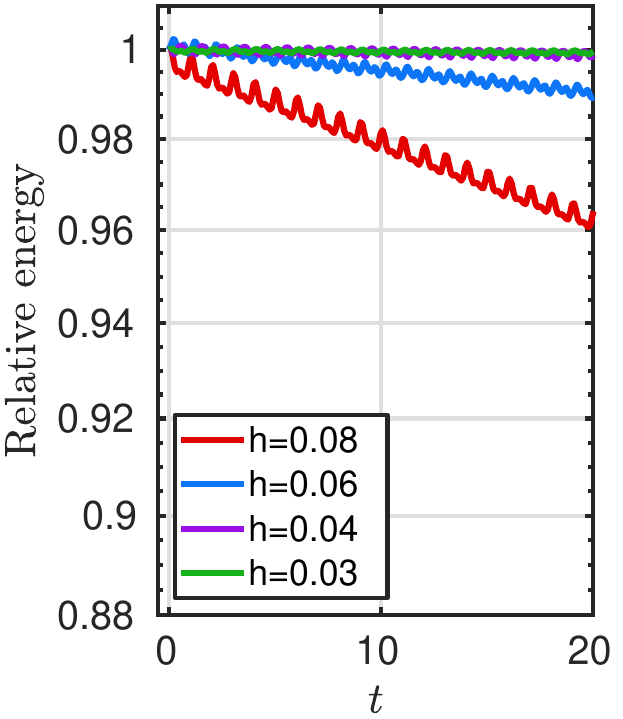} & 
    \hspace{-0.4cm}\includegraphics[width=0.3\linewidth]{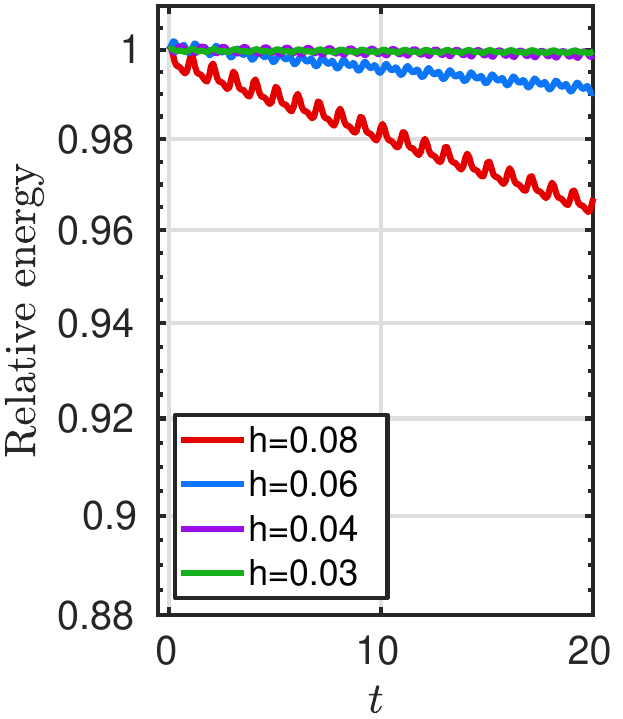} 
\end{tabular}
\caption{Energy as a function of time, relative to the energy of the initial condition (Kansa's RBF method), measured in the $\ell_2$-norm, for different choices of: the internodal distance $h$ (different lines in figures), the oversampling parameter $q$ (columns) and 
the monomial basis degree $p=5$, used to construct the global approximation. $\text{CFL}=0.6$ was used to determine the time step.}
\label{fig:experiments:conservation_energy:Kansa}
\end{figure}

\begin{figure}
    \centering
\begin{tabular}{lccc}
    \multicolumn{4}{c}{\hspace{0.8cm}\textbf{RBF-PUM: energy in time} \vspace{0.1cm}} \\
    & \hspace{0.6cm} $\mathbf{q=6}$ & \hspace{0.5cm}$\mathbf{q=9}$ & \hspace{0.5cm} $\mathbf{q=12}$\\
    \rotatebox{90}{\hspace{2.2cm}$\mathbf{p=2}$}  & \includegraphics[width=0.29\linewidth]{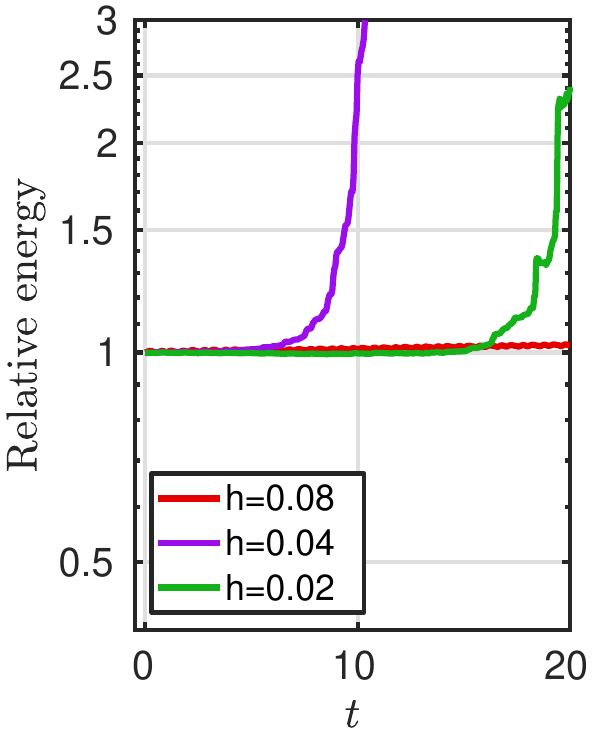} & 
    \hspace{-0.4cm}\includegraphics[width=0.3\linewidth]{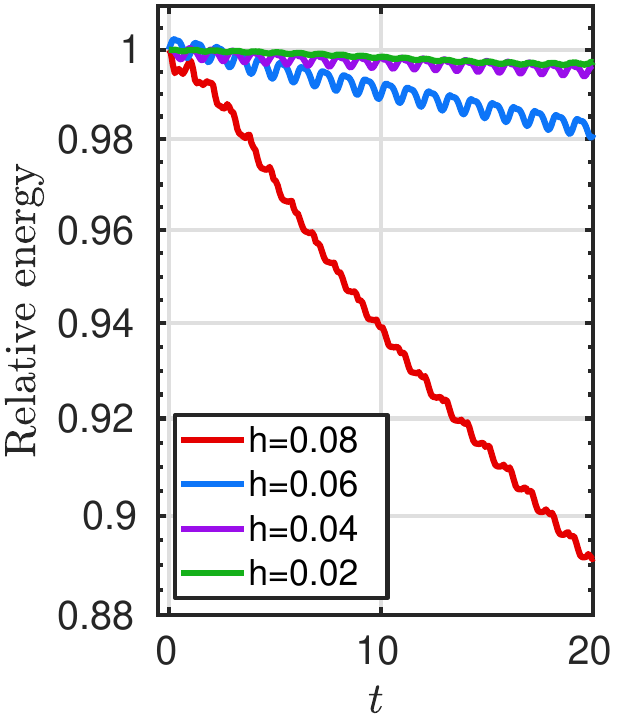} & 
    \hspace{-0.4cm}\includegraphics[width=0.3\linewidth]{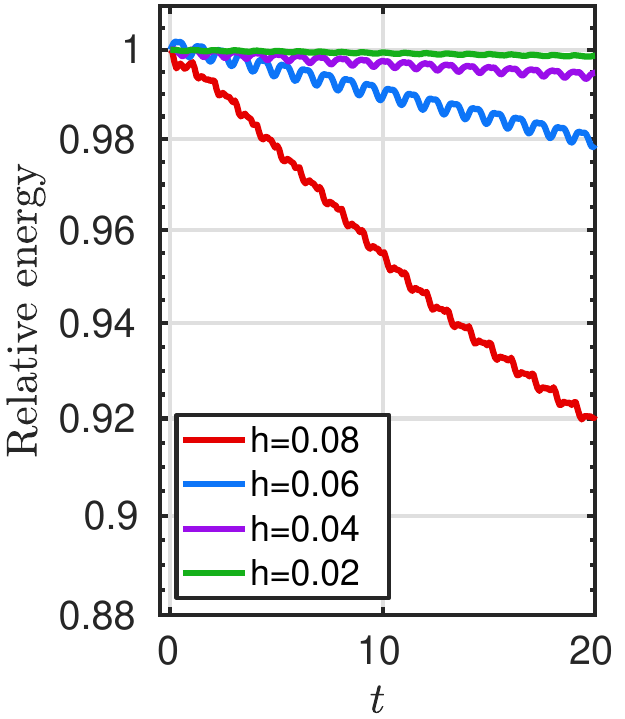} \\
    & \hspace{0.8cm} $\mathbf{q=6}$ & \hspace{0.5cm}$\mathbf{q=9}$ & \hspace{0.5cm} $\mathbf{q=12}$\\
    \rotatebox{90}{\hspace{2.2cm}$\mathbf{p=5}$}  & \includegraphics[width=0.3\linewidth]{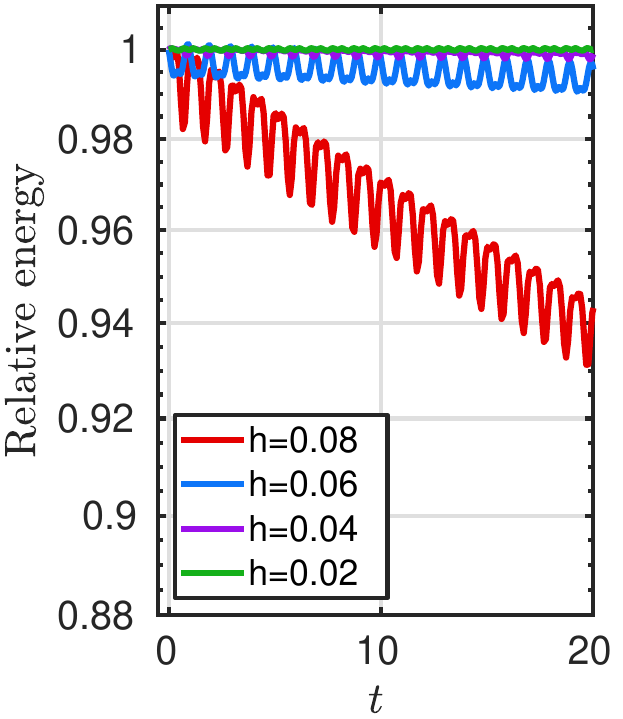} & 
    \hspace{-0.4cm}\includegraphics[width=0.3\linewidth]{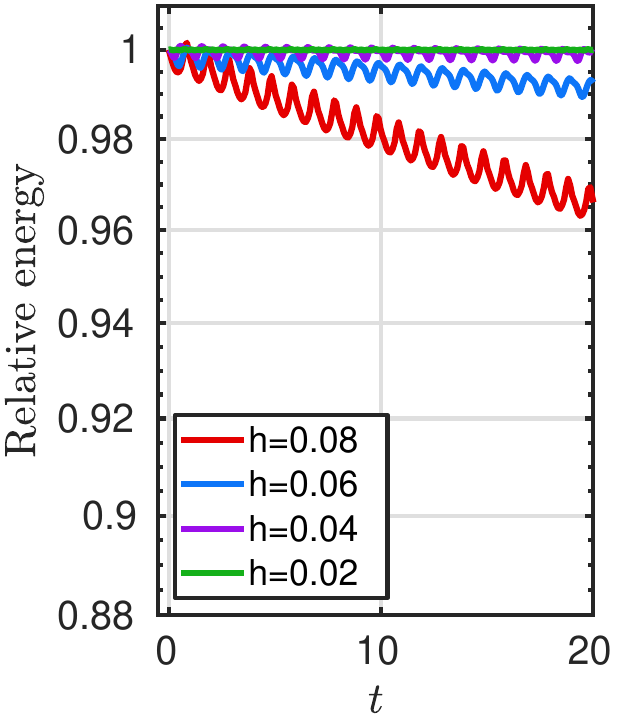} & 
    \hspace{-0.4cm}\includegraphics[width=0.3\linewidth]{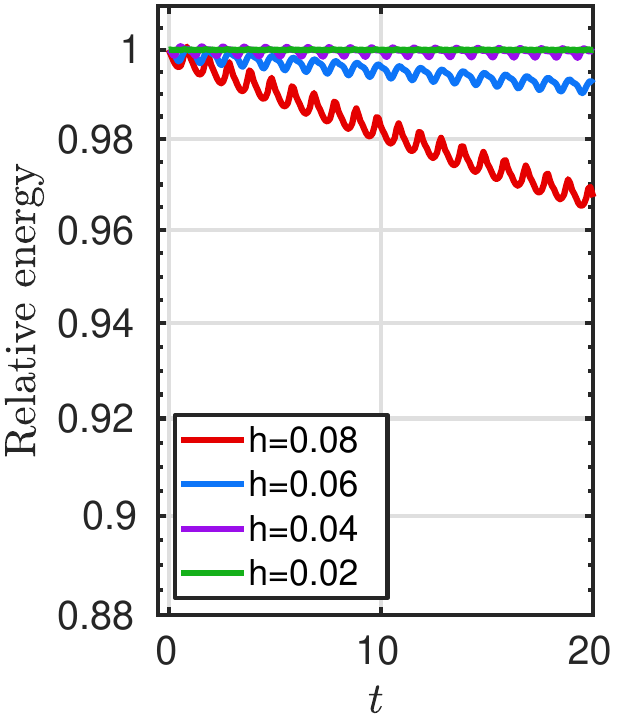} 
\end{tabular}
\caption{Energy as a function of time, relative to the energy of the initial condition (RBF-PUM), measured in the $\ell_2$-norm, for different choices of: the internodal distance $h$ (different lines in figures), the oversampling parameter $q$ (columns) and 
the monomial basis degrees $p$ (rows) which are 
used to construct the patch-based approximation. $\text{CFL}=0.6$ was used to determine the time step.}
\label{fig:experiments:conservation_energy:rbfpum}
\end{figure}

The result for Kansa's RBF method is displayed in Figure \ref{fig:experiments:conservation_energy:Kansa}, where we fixed $p=5$. 
Firstly, we observe that the method is 
stable for all considered parameters, as the energy is either constant or decaying. This is in agreement with the semi-discrete stability 
estimate \eqref{eq:stabilityKansa:finalEstimate}. Secondly, the energy decays when $h$ is large, and then gets asymptotically more 
constant as $h$ is decreased. Thirdly, the oversampling does not significantly improve the conservation of energy.
Results for RBF-PUM are collected in Figure \ref{fig:experiments:conservation_energy:rbfpum}. The results are very similar to the Kansa's method case, except that 
the method is unstable for some $h$ and $p$ when the oversampling is too small ($q=6$). 
The method does get stable when the oversampling is larger ($q=9$, $q=12$), which is in agreement 
with the stability estimate derived in Section \ref{sec:stabilityproperties:rbfpum}. 
We also observe that RBF-PUM is slightly more diffusive when $p=2$ and $h$ is small, compared with Kansa's RBF method.

\begin{figure}
    \centering
\begin{tabular}{lccc}
    \multicolumn{4}{c}{\hspace{0.8cm}\textbf{The unstabilized RBF-FD method: energy in time} \vspace{0.1cm}} \\
    & \hspace{0.6cm} $\mathbf{q=6}$ & \hspace{0.4cm}$\mathbf{q=9}$ & \hspace{0.4cm} $\mathbf{q=12}$\\
    \rotatebox{90}{\hspace{2.2cm}$\mathbf{p=2}$}  & \includegraphics[width=0.3\linewidth]{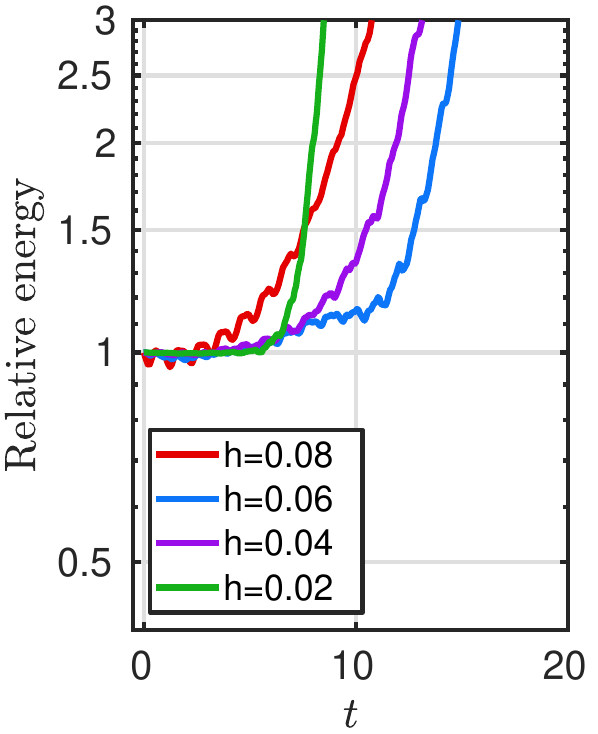} & 
    \hspace{-0.4cm}\includegraphics[width=0.3\linewidth]{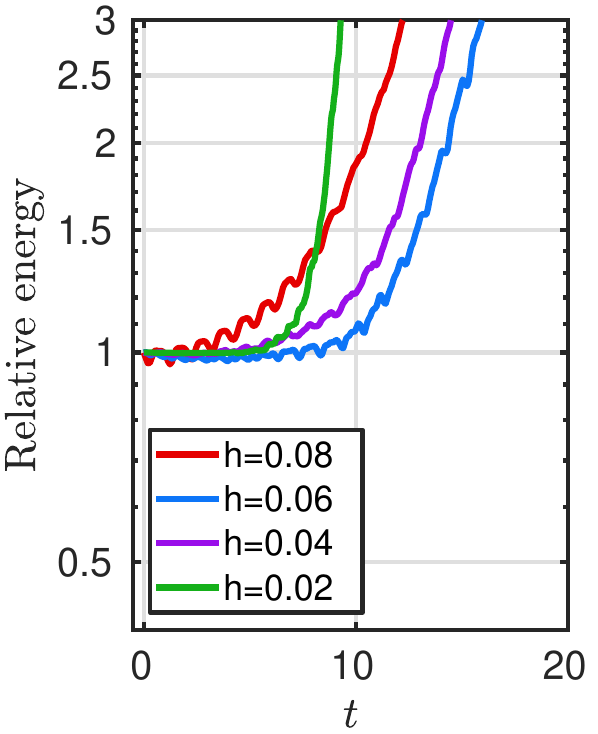} & 
    \hspace{-0.4cm}\includegraphics[width=0.3\linewidth]{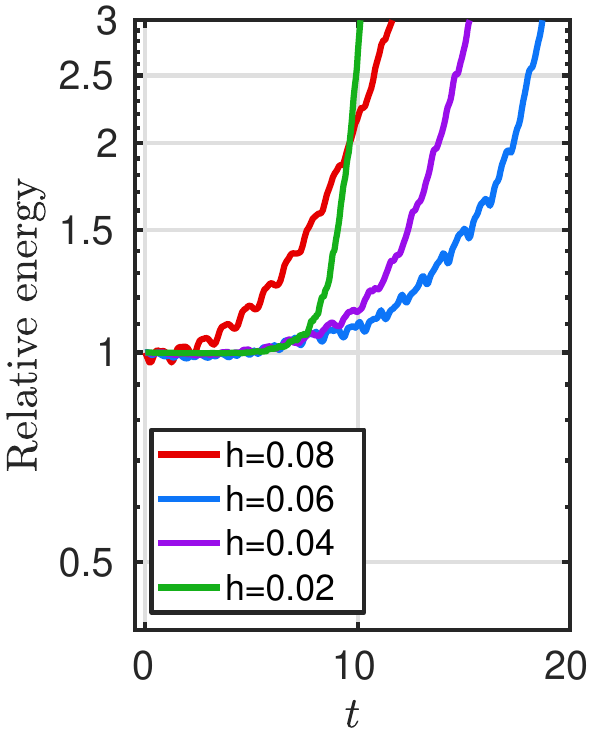} \\
    & \hspace{0.6cm} $\mathbf{q=6}$ & \hspace{0.5cm}$\mathbf{q=9}$ & \hspace{0.5cm} $\mathbf{q=12}$\\
    \rotatebox{90}{\hspace{2.2cm}$\mathbf{p=5}$}  & \includegraphics[width=0.3\linewidth]{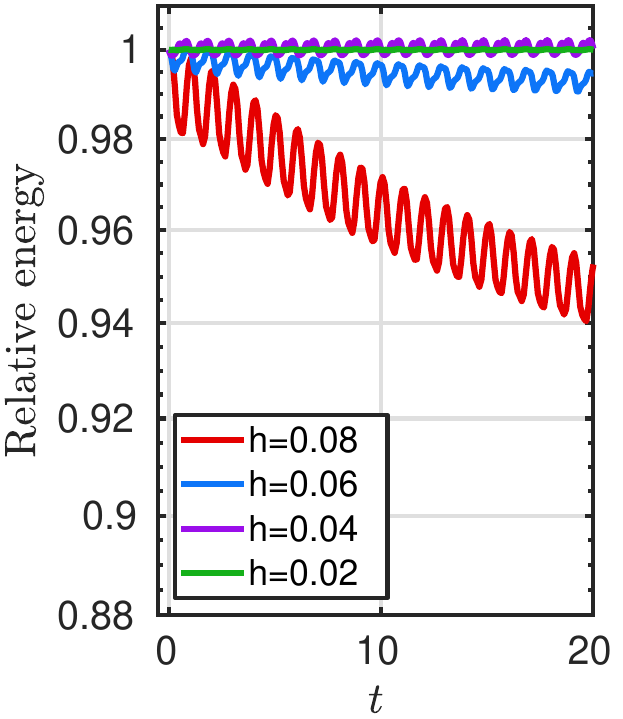} & 
    \hspace{-0.4cm}\includegraphics[width=0.3\linewidth]{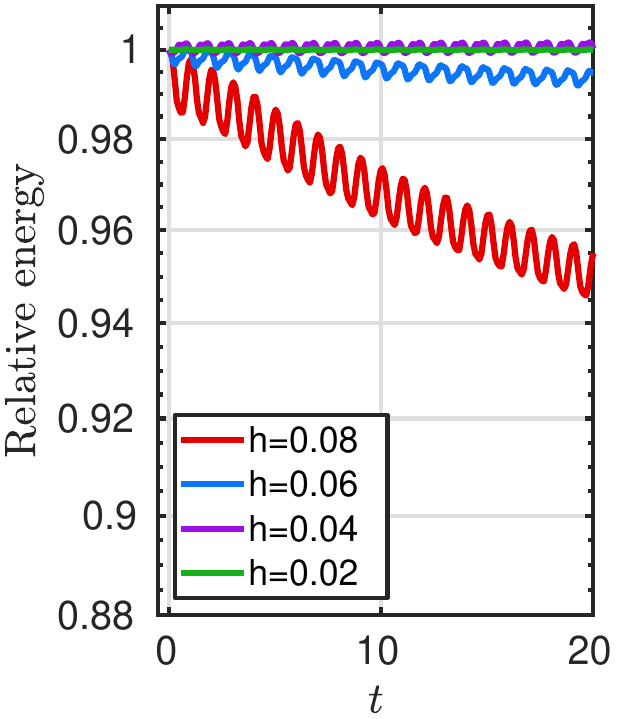} & 
    \hspace{-0.4cm}\includegraphics[width=0.3\linewidth]{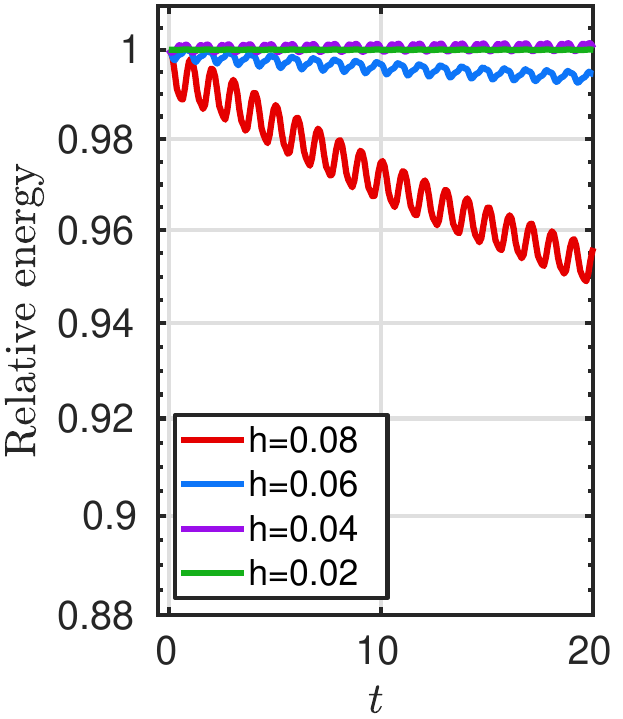} 
\end{tabular}
\caption{Energy as a function of time, relative to the energy of the initial condition (unstabilized RBF-FD), measured in the $\ell_2$-norm, for different choices of: the internodal distance $h$ (different lines in figures), the oversampling parameter $q$ (columns) and 
the monomial basis degrees $p$ (rows) which are 
used to construct the stencil-based approximation. $\text{CFL}=0.6$ was used to determine the time step.}
\label{fig:experiments:conservation_energy:rbffd}
\end{figure}

The unstabilized RBF-FD method is examined through Figure \ref{fig:experiments:conservation_energy:rbffd}. 
Here we observe that the method is unstable for all $h$ and $q$ when $p=2$. This is in agreement with the stability estimate 
\eqref{eq:stabilityrbffd:finalEstimate}, but also with the examination of the stability of the eigenvalue spectra in Figure 
\ref{fig:experiments:eigenvalues:rbffd}. Further, when $p=2$ we observe that the method does not blow up until $t=3$, 
at which point we made the results in Figure \ref{fig:experiments:error_href}. 
When $p=5$ the method is stable for all choices of the considered parameters, which is due 
to an increase in $n$ (compared with the case $p=2$) and smaller jumps as a consequence, 
as already pointed out in Section \ref{sec:experiments:eigenvalues} through Figure \ref{fig:experiments:eigenvalues:rbffd}. We emphasize that we can not guarantee 
the unstabilized RBF-FD method to be stable when $p$ is large, as the choice of $n$ at which the jumps become small enough can depend on many other parameters.

\begin{figure}
    \centering
\begin{tabular}{lccc}
    \multicolumn{4}{c}{\hspace{0.8cm}\textbf{The jump-stabilized RBF-FD method: energy in time} \vspace{0.1cm}} \\
    & \hspace{0.8cm} $\mathbf{q=6}$ & \hspace{0.5cm}$\mathbf{q=9}$ & \hspace{0.5cm} $\mathbf{q=12}$\\
    \rotatebox{90}{\hspace{2.2cm}$\mathbf{p=2}$}  & \includegraphics[width=0.3\linewidth]{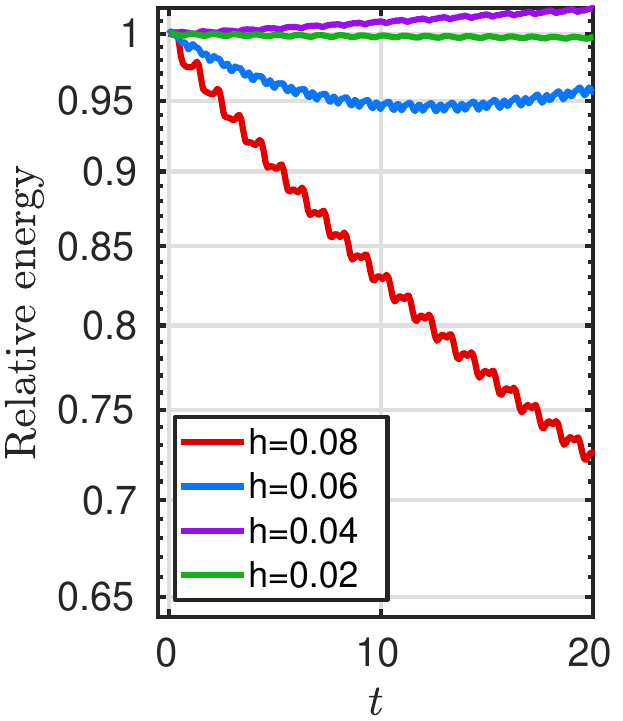} & 
    \hspace{-0.4cm}\includegraphics[width=0.3\linewidth]{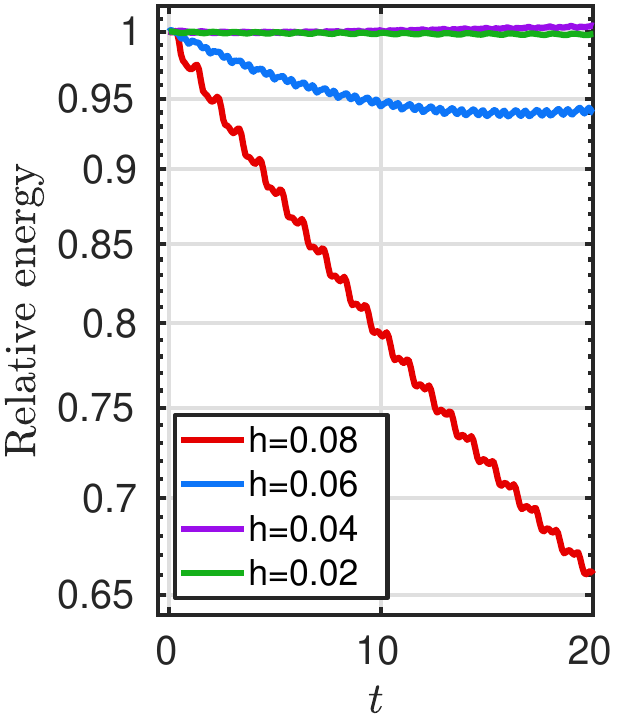} & 
    \hspace{-0.4cm}\includegraphics[width=0.3\linewidth]{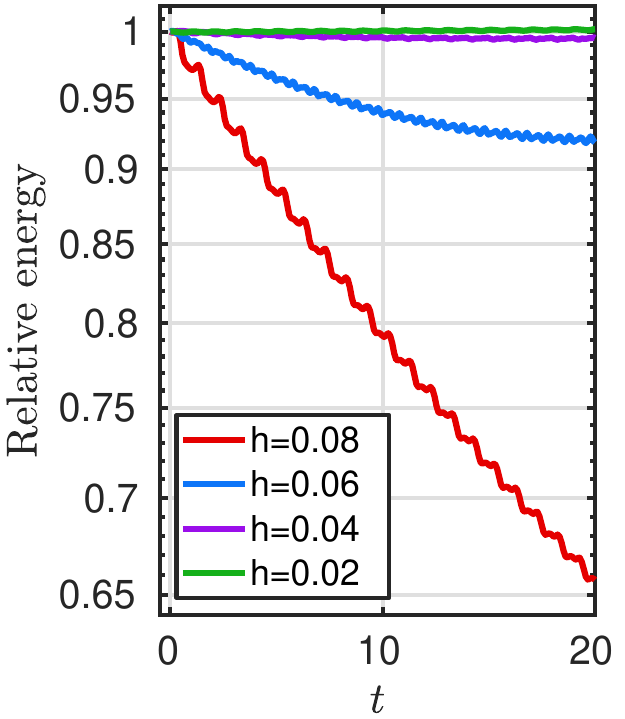} \\
    & \hspace{0.8cm} $\mathbf{q=6}$ & \hspace{0.5cm}$\mathbf{q=9}$ & \hspace{0.5cm} $\mathbf{q=12}$\\
    \rotatebox{90}{\hspace{2.2cm}$\mathbf{p=5}$}  & \includegraphics[width=0.3\linewidth]{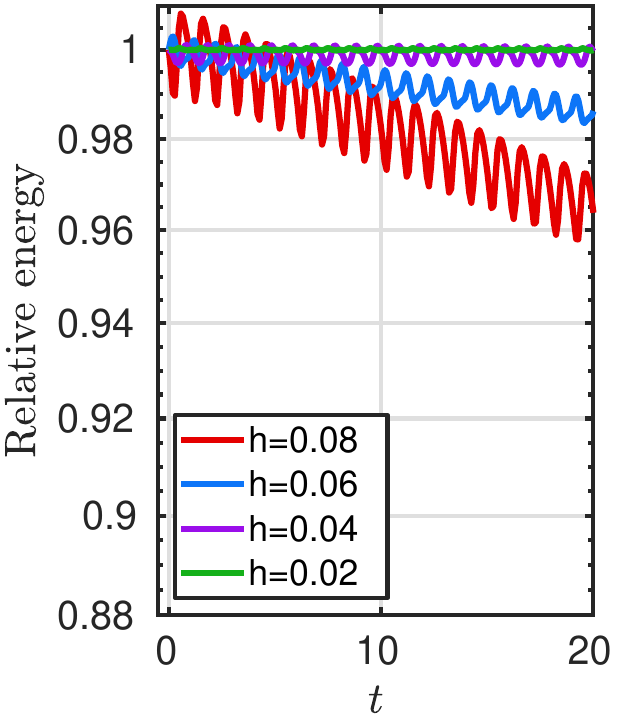} & 
    \hspace{-0.4cm}\includegraphics[width=0.3\linewidth]{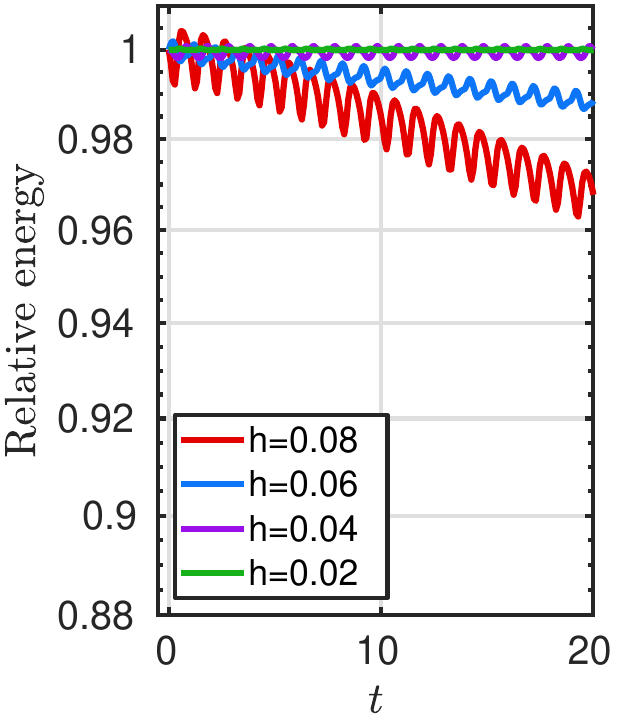} & 
    \hspace{-0.4cm}\includegraphics[width=0.3\linewidth]{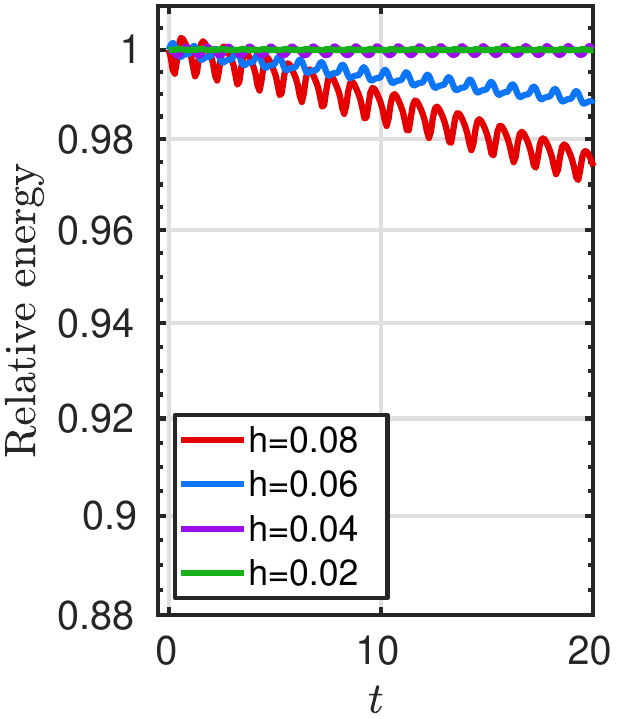} 
\end{tabular}
\caption{Energy as a function of time, relative to the energy of the initial condition (jump-stabilized RBF-FD), measured in the $\ell_2$-norm, for different choices of: the internodal distance $h$ (different lines in figures), the oversampling parameter $q$ (columns) and 
the monomial basis degrees $p$ (rows) which are 
used to construct the stencil-based approximation. $\text{CFL}=0.6$ was used to determine the time step.}
\label{fig:experiments:conservation_energy:rbffd_jumpstabilized}
\end{figure}

The results for the jump-stabilized RBF-FD method are given in Figure \ref{fig:experiments:conservation_energy:rbffd_jumpstabilized}. For all $p$, $q$ and $h$, 
the method is stable in time. When $p=2$, the method is diffusive for a large $h$. When $p=5$, the method is less 
diffusive for a large $h$. We observe that as $h$ is decreased, the conservation of energy is being asymptotically established. 

\subsection{Maximal CFL numbers for Kansa's RBF method, RBF-PUM and the RBF-FD method}
For all three methods, we examined the maximal CFL numbers that still allow stable time stepping when using the explicit classical Runge-Kutta 4 method 
to advance the solution in time. We are in particular interested in how the oversampling influences the CFL numbers and whether the jump-stabilization 
for the RBF-FD method has a negative effect on the maximal time step. 
The relation between $\Delta t_n$ and CFL is given in \eqref{eq:experiments:timestep}.
We fix the mean internodal distance in the perturbed DistMesh point set X to $h=0.04$ and compute the maximal CFL numbers as a function of the 
oversampling parameters $q=5,6,...,14$ and the monomial basis degrees $p=2,3,4,5,6$.
It turns out that for all the methods the CFL number has to be taken smaller 
as the oversampling parameter increases, however, the decay in the CFL number is small. 
The most significant difference between the smallest and the largest 
CFL numbers per a fixed $p$ is only $0.07 \%$, measured relatively to the largest CFL number per a fixed $p$. 
Thus, an increase in the oversampling does not catastrophically impact the size of $\Delta t$.

\section{Final remarks}
\label{sec:finalremarks}
In this paper we investigated the stability properties of Kansa's RBF method, RBF-PUM and the RBF-FD method, when applied 
to a time-dependent hyperbolic problem. Our analysis was made for a linear advection problem with a divergence free velocity field. 

We derived that the $\ell_2$-stability in time for Kansa's RBF method and RBF-PUM, is for a fixed $h$, controllable by $h_y$ which is inversely proportional to the oversampling parameter $q$. 
When collocation is used instead of oversampling, then $q=1$, and the stability control through $q$ is lost.
The $\ell_2$-stability in time of the RBF-FD method also depends on $q$ (controllable), 
but also on the jump terms across the interfaces of the Voronoi cells centered at each $x_i \in X$.
As we demonstrated numerically, a possibility to control the size of the jumps is to increase the stencil size $n$, however, this is a 
computationally expensive 
approach since an increase to $n$ gives an increase in the density of the discretization matrices. 
We provided a special penalty term that controls the jump term for any $n$, and in this way stabilizes the method.

Numerical experiments confirmed that by increasing the oversampling, the stability of Kansa's RBF method and RBF-PUM was being 
manifested, since we could find a $\text{CFL}$ number such that the eigenvalue spectra were fully contained in the classical explicit Runge-Kutta 4 stability region. 
The experiments also confirmed that 
the RBF-FD method can not be stabilized only by increasing the oversampling, unless we also decrease the spurious jump term at the same time.

\red{We considered a fairly simple case of linear advective flow with divergence-free conditions. As this assumption makes the flux skew-symmetric, the resulting analysis was easy to obtain. However, the above analysis can be extended to more general fluxes, including nonlinear ones. For example, in \cite{Nordstrom_2006} the authors demonstrate 
a splitting approach in which the divergence of the linear flux is written as a linear combination of symmetric and antisymmetric parts. Then, by carefully selecting the coefficients in this linear combination, the authors obtain $\ell_2$-stability.
For more general nonlinear fluxes, we refer the reader to \cite{Tadmor_1984}, where it was shown that the skew-symmetric form of many nonlinear fluxes $\bm f(u)$ can be constructed if the corresponding entropy function and entropy flux of the PDE are homogeneous functions. 
Once the skew-symmetric form of the flux is obtained, the frozen coefficients approach \cite{Mishra_2010} can be applied to make the corresponding flux terms linear and obtain stability estimates. }

Overall, we established a theoretical framework and illustrated some of the practical implications, which can, in the future, 
be used for further improvements of the three considered RBF methods.

\appendix
\section{The integration error induced by the oversampling, for piecewise continuous functions}
\label{sec:appendix:integration}
\noindent
In this section we derive an estimate for the integration error 
\begin{equation}
\frac{|\Omega|}{M} \sum_{i=1}^M f(y_i) - \int_\Omega f\, d\Omega, 
\end{equation}
where $y_i \in \Omega \subset \mathbb{R}^2$ and $|\Omega| = \int_{\Omega}1\, d\Omega$.
%
The finite-dimensional piecewise smooth integrand $f \in L_2(\Omega)$. At the same time we also have that $f \in \mathcal C^k(K_i)$, $k>0$, $i=1,..,N$, 
where $\cup_{i=1}^N K_i = \Omega$ is the set of 
Voronoi regions, where each $K_i$ is as defined in \eqref{eq:voronoi_region_Ki}. In the estimate, we use regions that intersect the domain boundary. 
Therefore we introduce the modified integrand  $\bar f = f g$, where $g=1$ on $\Omega$ and $g=0$ elsewhere.
%
We limit the derivation to the case where the point set $Y$ of size $M$ has a 
Cartesian layout with grid size and internodal distance $h_y$, subjected to the domain $\Omega$. We let $y_{\Gamma(i,j)} \in Y$, $j=1,..,q_i$ be the points that sample the Voronoi region $K_i$, 
where $\Gamma$ is an operator relating the local index of an evaluation point contained in $K_i$, to a global index in $Y$. For a quasi uniform node set $X$, $q_i=\alpha_iq$ with $\alpha_i\approx1$. 
We bound the integration error using boxes $B_{ij}$ with side $h_y$, centered at $y_{\Gamma(i,j)}$. Out of the $q_i$ boxes there are $n_i$ boxes such that $B_{ij}\subset K_i$, while $q_i-n_i$ 
extend outside $K_i$. Due to the area to boundary relation, the box numbers can be bounded by $n_i\leq c_Aq$ and $(q_i-n_i) \leq c_B\sqrt{q}$. 
Furthermore, we need to add some auxiliary boxes $\tilde B_{ij}$ that do not have a corresponding center point in $Y$, but are needed to cover all of $\Omega$. 
In the same way, the number $n_\Omega$ of these boxes can be bounded by 
$n_\Omega\leq c_\Omega \sqrt{N}\sqrt{q}$, since the number of interior boxes $B_{ij}$ is $M=Nq$.

The numerical integration over one box is written as 
$I_h(B_{ij})=\frac{|\Omega|}{M} \bar f(y_{\Gamma(i,j)})=|B_{ij}|\bar{f}(y_{\Gamma(i,j)})=\int_{B_{ij}}\bar{f}(y_{\Gamma(i,j)})\,dy$. 
Furthermore, we note that due to the continuity of $f$ on $K_i$, the mean value theorem leads to $\bar{f}(y)= \bar f(y_{\Gamma_{ij}})+ (y-y_{\Gamma_{ij}}) \cdot \nabla \bar f(\xi)$, where $\xi=\theta y +(1-\theta)y_{\Gamma_{ij}}$, $0\leq\theta\leq 1$, $y\in K_i$.
The integration error over $B_{ij}$ then becomes:
\begin{equation}
    \begin{aligned}
  I(B_{ij})-I_h(B_{ij})=\int_{B_{ij}} \bar f(y)-\bar f(y_{\Gamma_{ij}})\,dy
  &= \int_{B_{ij} \cap K_{i}}(y-y_{\Gamma_{ij}}) \cdot \nabla \bar f(\xi(y))\,dy\nonumber\\
  & + \int_{B_{ij}\setminus (B_{ij} \cap K_i)}\bar f(y)-\bar f(y_{\Gamma_{ij}})\,dy. 
    \end{aligned}
\end{equation}
The integration error over the whole domain can be estimated as:
\begin{equation}
    \begin{aligned}
  |I(\Omega)-I_h(\Omega)| & =  \left|\sum_{i=1}^{n_\Omega}\int_{\tilde{B}_{ij}}\bar f(y)\,dy + \sum_{i=1}^N\sum_{j=1}^{q_i}I(B_{ij})-I_h(B_{ij})\right|\\
  & \leq  \sum_{i=1}^{n_\Omega}\left|\int_{\tilde{B}_{ij}} \bar f(y)\,dy\right| + \sum_{i=1}^N\sum_{j=1}^{q_i} \left|\int_{B_{ij} \cap K_{i}}(y-y_{\Gamma_{ij}}) \cdot \nabla \bar f(\xi(y))\,dy\right| \\
  &\quad + \sum_{i=1}^N\sum_{j=1}^{q_i-n_i}\left|\int_{B_{ij}\setminus (B_{ij} \cap K_i)}\bar f(y)-\bar f(y_{\Gamma_{ij}})\,dy\right| \\
  &\leq c_\Omega\sqrt{N}\sqrt{q} \max_{ij} |\tilde{B}_{ij}| \max_{ij} \|f\|_{L_\infty(\tilde B_{ij})} + \\
  & \quad + Nq \max_i\alpha_i \max_{ij} |B_{ij}\cap K_i|\frac{h_y}{\sqrt{2}} \max_{ij} \|\nabla f\|_{L_\infty(B_{ij}\cap K_i)}\\
  & \quad + c_B N \sqrt{q} \max_{ij} |B_{ij}\setminus (B_{ij} \cap K_i)| 2 \max_{ij} \|f\|_{L_\infty(B_{ij}\setminus (B_{ij} \cap K_i))}\\
  &\leq (\sqrt{|\Omega|} c_\Omega   + 2 |\Omega| c_B h^{-1}) h_y\, \max_i \|f\|_{L_\infty(K_i)}  + \\
  & \quad + |\Omega| \frac{\max_i\alpha_i}{\sqrt{2}}h_y \max_{i} \|\nabla f\|_{L_\infty(K_i)},
\end{aligned}
\end{equation}
where we in the last step used $|B_{ij}|=h_y^2$, $|B_{ij}\cap K_i|\leq |B_{ij}|$,  $|B_{ij}\setminus (B_{ij} \cap K_i)|\leq |B_{ij}|$, $\sqrt{N}\sqrt{q}=\sqrt{|\Omega|}\, h_y^{-1}$, and $1/\sqrt{q}=\frac{h_y}{h}$. 
We finish the estimate by using an inverse inequality $\max_{i} \|\nabla f\|_{L_\infty(K_{i})} \leq C_I h^{-1} \max_{i} \| f \|_{L_\infty(K_{i})}$ from \cite{BrennerBook},  where $C_I$ is independent of $h$ to arrive to:
\begin{equation}
    \label{appendix:integral:Omega}
    \begin{aligned}
    |I(\Omega)-I_h(\Omega)| &\leq (\sqrt{|\Omega|} c_\Omega   + 2 |\Omega| c_B h^{-1} + |\Omega| \frac{\max_i\alpha_i}{\sqrt{2}} C_I h^{-1} ) h_y\, \max_i \|f\|_{L_\infty(K_i)} \\
        &\leq C_{\int}\, h^{-1} h_y \max_i \|f\|_{L_\infty(K_i)},
    \end{aligned}
\end{equation}
where $C_{\int} = \sqrt{|\Omega|} c_\Omega   + |\Omega| (2c_B + \frac{\max_i\alpha_i}{\sqrt{2}}C_I)$. 
We observe that for a fixed $h$, the integration error induced by oversampling asymptotically approches $0$ as $h_y \to 0$.

\section*{Acknowledgments}
We thank Víctor Bayona from Universidad Carlos III de Madrid and Katharina Kormann from Ruhr University Bochum for fruitful discussions.


\bibliography{refs}

\end{document}